\documentclass[11pt,a4paper,leqno,twoside]{amsart}
\usepackage{latexsym,amssymb,amsmath}
\usepackage{epsfig}
\input xy
\xyoption{all}
\usepackage[english]{babel}

\usepackage{amsthm}
\usepackage{amscd}
\usepackage{tabularx}

\def\o{\omega}

\def\cS{{\mathcal S}}

\def\CC{\mathbb C}
\def\RR{\mathbb R}
\def\HH{\mathbb H}
\def\AA{{\mathbb A}}
\def\BB{{\mathbb B}}
\def\OO{\mathbb O}
\def\ZZ{\mathbb Z}

\def\PP{\mathbb P}

\def\FF{\mathbb F}

\def\fh{{\mathfrak h}}

\def\fsl{{\mathfrak {sl}}}

\def\fspin{{\mathfrak {spin}}}
\def\fso{{\mathfrak {so}}}
\def\fe{{\mathfrak e}}

\def\ff{{\mathfrak f}}

\def\fg{{\mathfrak g}}

\def\fz{{\mathfrak z}}
\def\ft{{\mathfrak t}}

\def\a{\alpha}

\def\b{\beta}
\def\g{\gamma}
\def\s{\sigma}

\def\d{\delta}
\def\th{\theta}

\def\ot{{\mathord{\,\otimes }\,}}
\def\op{{\mathord{\,\oplus }\,}}

\def\ra{{\mathord{\;\rightarrow\;}}}

\newtheorem{theo}{Theorem}

\newtheorem{lemm}[theo]{Lemma}
\newtheorem{prop}[theo]{Proposition}

\begin{document}

\title{Configurations of lines \linebreak 
and models of Lie algebras}
\author{L. Manivel}
\maketitle

\begin{abstract}
The automorphism groups of the $27$ lines on the smooth 
cubic surface or the $28$ bitangents to the general quartic plane 
curve are well-known to be closely related to the Weyl groups 
of $E_6$ and $E_7$. We show how classical subconfigurations
of lines, such as double-sixes, triple systems or Steiner sets, 
are easily constructed from certain models of the exceptional 
Lie algebras. For $\fe_7$ and $\fe_8$ we are lead to beautiful 
models graded over the octonions, which display these algebras as 
plane projective geometries of subalgebras. We also interprete
the group of the bitangents as a group of transformations of the 
triangles in the Fano plane, and show how this is related to an 
interpretation of the isomorphism $PSL(3,\FF_2)\simeq PSL(2,\FF_7)$
in terms of harmonic cubes.
\end{abstract}

\section{Introduction}

Such classical configurations of lines as the $27$ lines 
on a complex cubic surface or the $28$ bitangents to a smooth 
quartic plane curve have been extensively studied in the 19th 
century (see e.g. \cite{hen}). Their automorphism groups were known, but 
only at the beginning of the 20th century were their close 
connections with the Weyl groups of the root systems $E_6$ and 
$E_7$, recognized, in particular through the relation with Del 
Pezzo surfaces of degree three and two, respectively \cite{cox,dem}. 
Del Pezzo surfaces of degree one provide a similar identification of 
the diameters of the root system $E_8$, with the $120$ tritangent 
planes to a canonical space curve of genus $4$. 

Can we go beyond the Weyl groups and find a connection with 
the Lie groups themselves? The $27$ lines on the cubic surface
are in natural correspondence with the weights of the minimal 
representation of $E_6$, from which the Lie group can be 
recovered as the stabilizer of a cubic form that already appears in Elie
Cartan's thesis; Faulkner showed
how to define this form in terms of the $45$ tritangent planes. 
A similar phenomenon holds for the $28$ bitangents to the quartic 
plane curve, which can be put in correspondence with pairs of opposite 
weights of the minimal representation of $E_7$. In both cases the 
connection between the Lie group and its Weyl group is particularly 
close because of the existence of a minuscule representation. 
For $E_8$ the minimal representation is the adjoint one and is no
longer minuscule. 

The first aim of this paper is to use these connections with the 
Lie groups, or rather the Lie algebras $\fe_6$, $\fe_7$, $\fe_8$,
to shed a new light on the work of the classical geometers on 
line configurations. Our main idea is that each time we consider 
a semisimple Lie subalgebra, the restriction of the minimal 
representation branches into a direct sum of subrepresentations,
and consequently the weights split into special subsets forming 
interesting subconfigurations. In the case of $\fe_6$ and the $27$
lines we get the following correspondence:
\smallskip

\begin{center}
\begin{tabular}{ccc}
{\it Subalgebra} & \hspace*{1cm} \qquad & {\it Subconfiguration} \\
$\fspin_{10}$ & & Line \\
$\fspin_{8}$ & & Tritangent plane \\
$\fsl_2\times\fsl_6$ & & Double-six \\
$\fsl_3$ & & Steiner set \\
$\fsl_3\times\fsl_3\times\fsl_3$ & & Triple system 
\end{tabular}
\end{center}

\smallskip
In the case of $\fe_7$ and the $28$ bitangents the notion of 
Steiner complexes of bitangents make a natural appearance. They
are special sets of $12$ bitangents which can be put in correspondence 
with positive roots, and also with points of a $5$-dimensional
projective space over $\FF_2$. This space is endowed with a natural
symplectic form, and indeed the Weyl group of $E_7$ is closely connected
with the finite symplectic group $Sp(6,\FF_2)$. This leads to a very
interesting finite symplectic geometry whose lines are known to 
correspond to the so called {\it syzygetic triads} of Steiner sets.
We prove that planes in this geometry are in correspondence with what 
we call {\it Fano heptads} of bitangents. The upshot is a finite 
geometry modeling the symplectic geometries related to the third line 
of Freudenthal Magic Square, whose last term is precisely $\fe_7$
\cite{LMfreud}. 

Very interestingly, this leads to a beautiful model of $\fe_7$
and its minimal representation which, rather unexpectedly, turns out 
to be closely related with the Fano plane and the octonionic
multiplication. Indeed, recall that $\OO$, the Cayley algebra of 
octonions, can 
be defined as the eight-dimensional algebra with a basis $e_0=1,
e_1,\ldots,e_7$, with multiplication rule encoded in an oriented 
Fano plane. 

\begin{center}
\setlength{\unitlength}{4mm}
\begin{picture}(30,11)(-15,-1)
\put(-5.5,0){ \resizebox{!}{1.5in}{\includegraphics{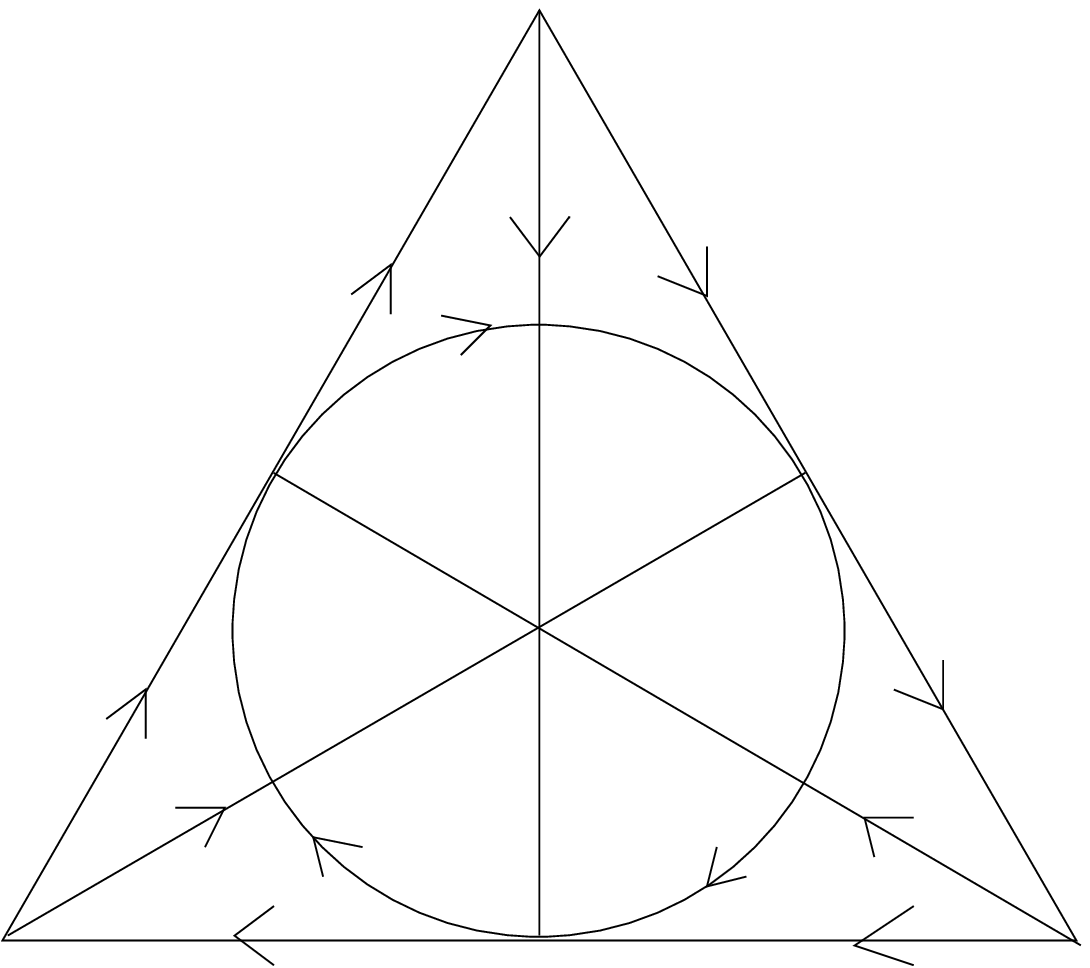}}}
\put(-3.6,4.6){$e_2$}\put(3.3,4.6){$e_6$}
\put(-.3,-.8){$e_4$}\put(-6.1,-.8){$e_1$}
\put(5.5,-.8){$e_5$}\put(-.2,9.7){$e_3$}\put(.9,3){$e_7$}
\end{picture} 
\end{center}

\centerline{{\it Figure 1. Octonionic multiplication}}

\bigskip
This means that $e_ie_j=\pm e_k$ if $i,j,k$ are three distinct points
on one of the projective lines in this plane, with a plus sign if and
only if $(ijk)$ gives the cyclic orientation fixed on the line.

\smallskip 
We define an $\OO$-grading on a Lie algebra $\fg$ to be a 
decomposition 
$$\fg =\fh_0 e_0\op\bigoplus_{1\le i\le 7}\fh_i e_i$$
such that $[\fh_i,\fh_j]\subset\fh_k$ if $e_ie_j=\pm e_k$. 
In particular $\fh_0$ is a subalgebra and each $\fh_i$ is an
$\fh_0$-module. More is true: for any point $i$ and any line 
$\ell$ in the Fano plane, the direct sums 
$$\fg_i=\fh_0\op\fh_i, \qquad 
\fg_{\ell} =\fh_0 e_0\op\bigoplus_{j\in\ell}\fh_j e_j$$
are subalgebras of $\fg$, so that we really have a configuration 
of Lie algebras defined by a plane projective geometry. 

Our discussion of Fano heptads lead us to discover that $\fe_7$
has a natural structure of an $\OO$-graded algebra, compatible 
with its action on the minimal representation $V$. Indeed, attach
to each line $\ell$ of the Fano plane a two dimensional
vector space $A_{\ell}$. Then we 
can describe $\fe_7$ and $V$ as follows:

\begin{eqnarray}
\nonumber 
\fe_7 &= &\times_{\ell}\fsl(A_{\ell})e_0\;\op\bigoplus_{1\le i\le 7}
(\otimes_{i\notin\ell}A_{\ell}) \, e_i, \\
\nonumber 
V &= &\bigoplus_{1\le j\le 7}(\otimes_{j\in\ell}A_{\ell}) \, e_j.
\end{eqnarray}

Going a little deeper in the Lie algebra structure, we will
discover a natural connection with the multiplication table
of the Cayley algebra. This leads to an amusing interpretation 
of the isomorphism $PSL(3,\FF_2)\simeq PSL(2,\FF_7)$ in terms 
of {\it harmonic cubes}, and a permutation representation 
of the group of the bitangents on the triangles of the Fano plane. 
  
A similar description of $\fe_8$ also exists, and the 
biggest two exceptional Lie algebras appear as plane projective
geometries whose points are copies of $\fso_8$ and $\fso_8\times\fso_8$,
and whose lines are copies of $\fso_{12}$ and $\fso_{16}$,
respectively. Moreover, this octonionic model of $\fe_8$ makes 
obvious the existence of the multiplicative orthogonal decomposition
that was a key ingredient in Thompson's construction 
of the sporadic simple group denoted $Th$ or $F_3$ (see 
\cite{kos}, Chapters 3 and 13). It would 
certainly be interesting to use this octonionic model, suitably adapted, 
to construct forms of $\fe_8$ over arbitrary fields. 
  
\medskip Classically, two unifying perspectives on the line
configurations we are interested in have been particularly successful. 
We briefly discuss the connection with our present approach. 

\smallskip\noindent {\bf Theta characteristics.} Bitangents to the 
plane quartic curve (a canonical curve of genus $g=3$), as well as
tritangent planes to the canonical curve of genus $g=4$,  
can be interpreted as odd theta-characteristics. Since the
theta-characteristics can be seen as points of an affine space
over the half-periods of the curve, this leads to an
interpretation in terms of finite symplectic geometries in dimension
$2g$ over the field $\FF_2$. This was developped in great detail by 
the classical geometers, in particular by 
Coble (\cite{coble}, Chapter II). For example the theta-characteristics 
can be understood as the quadrics whose associated polarity is the 
natural symplectic form. Isotropic linear spaces also have natural 
geometric interpretations. For $g=3$ one recovers the connection
that we already mentionned between $W(E_7)$ and $Sp(6,\FF_2)$. For 
$g=4$, the Weyl group $W(E_8)$ is the automorphism group of the
lines in the Del Pezzo surface of degree one, whose canonical 
model is a double covering of a quadratic cone, branched along 
a canonical sextic curve. As noticed by Schottky, there is a unique
even theta-characteristic vanishing at the vertex of the cone, which
explains why the automorphism group of the tritangent planes is 
an orthogonal group $O(8,\FF_2)^+$ rather than a symplectic group.

\smallskip\noindent {\bf Semi-regular polytopes.}
Gosset seems to have been the first, at the very beginning 
of the 20th century, to understand that the lines on the cubic 
surface can be interpreted as the vertices of a polytope, whose
symmetry group is precisely the automorphism group of the 
configuration. Coxeter extended this observation to the 
28 bitangents, and Todd to the 120 tritangent planes. Du Val
and Coxeter provided systematic ways to construct the polytopes, 
which are denoted $n_{21}$ for $n=2,3,4$ and live in $n+4$ 
dimensions \cite{duval,coxe,cox}. 
They have the characteristic property of being semi-regular, which 
means that the automorphism group acts transitively on the vertices, 
and the faces are regular polytopes. In terms of Lie theory they are
best understood as the polytopes in the weight lattices of the
exceptional simple Lie algebras $\fe_{n+4}$, whose vertices are 
the weights of the minimal representations. Coxeter investigated 
in great detail their semi-regular sub-polytopes \cite{coxe}. 
Algebraically, this amounts to identifying certain Lie 
subalgebras of the $\fe_{n+4}$. But Coxeter does not describe how
the full polytopes are organized around these special 
sub-polytopes. In a sense this is what we will be doing in this 
paper, with the nice conclusion that it leads to a very natural, 
unified and easy-going description of (at least part of) the 
classical combinatorics of the line configurations, as well as
new insights in the fascinating structure of the exceptional 
Lie algebras. 

\medskip\noindent {\it Acknowledgements:} I thank I. Dolgachev and
P.E. Chaput for their useful comments.

\section{Models of the exceptional Lie algebras}

\subsection*{The Reye configuration and triality}

A classical elementary configuration of lines is the Reye 
configuration below, obtained from a cube in a three 
dimensional projective space (see \cite{hilb} and \cite{dolg2}). 
This configuration can be understood as a central projection of the 
$24$-cell, one of the regular polytopes in four dimensions.
The vertices of this polytope are given by the roots of the 
root system $D_4$ (we use \cite{bou} as a general reference 
on root systems). 

\begin{center}
\setlength{\unitlength}{4mm}
\begin{picture}(30,17)(-10,-1)
\put(-5.5,0){ \resizebox{!}{2.5in}{\includegraphics{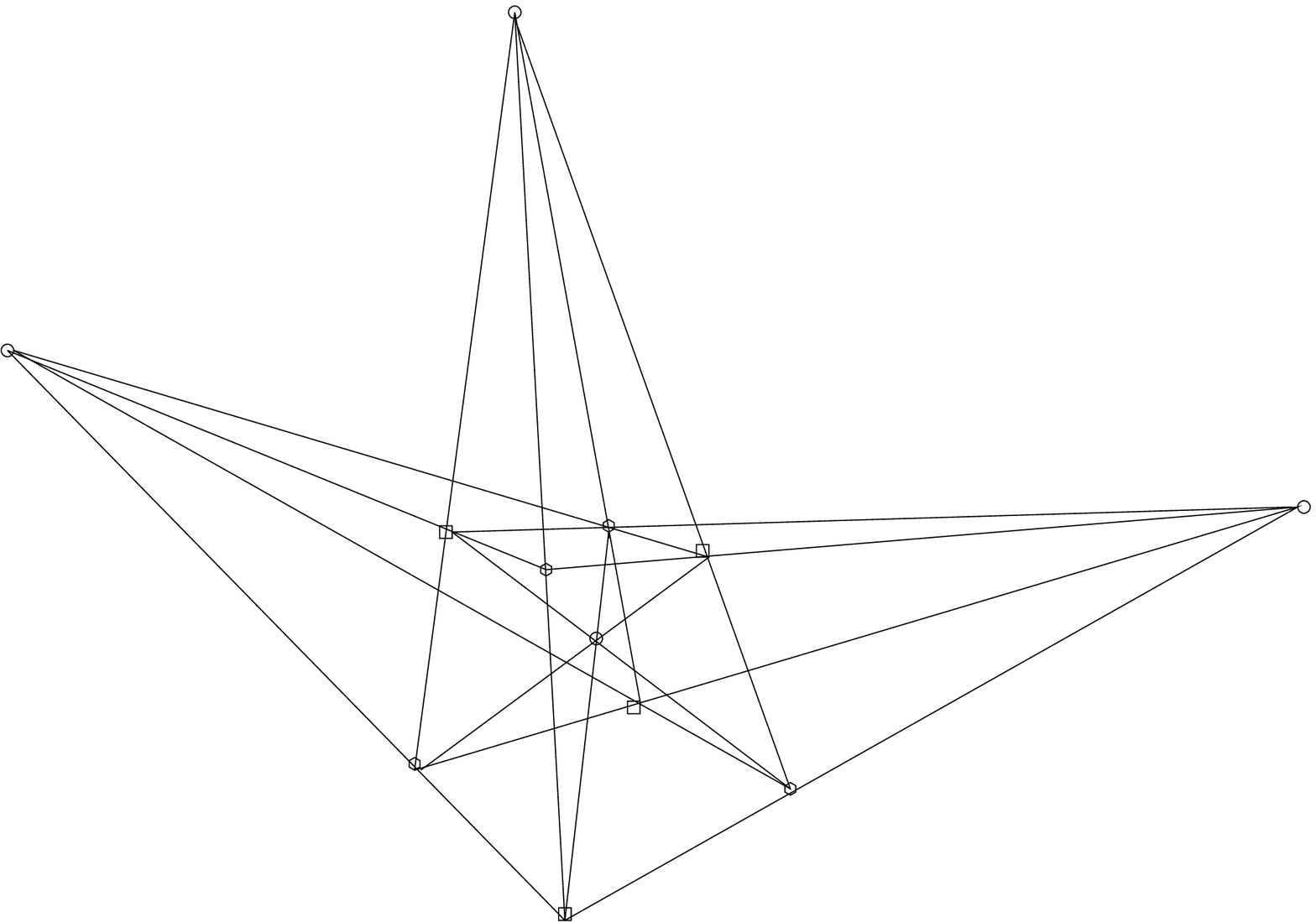}}}
\end{picture} 
\end{center}
\centerline{{\it Figure 2. The Reye configuration}}

\medskip
As one can easily see on Figure 2, there is a unique way to partition 
the points of the Reye configuration into three types, in such 
a way that each line contains exactly one point of each type. 
This decomposes the $24$-cell into three $16$-cells given by the 
vertices of three hypercubes. Each of these defines a root subsystem
of $D_4$ of type $A_1^4$. 

Restricting the adjoint representation 
of $\fspin_8$ to the corresponding subalgebra, a product of four 
copies of $\fsl_2$, we obtain the {\it four-ality model} \cite{LMpop}
$$\fspin_8=\fsl(A_1)\times\fsl(A_2)\times\fsl(A_3)\times\fsl(A_4)
\op (A_1\ot A_2\ot A_3\ot A_4),$$  
whose existence is indicated by the shape of the affine Dynkin diagram 
$\tilde{D}_4$. (Here $A_1, A_2, A_3, A_4$ are two dimensional 
complex vector spaces. Note that the construction works on the reals
to give the split form $\fso_{4,4}$.) 
Four-ality reduces to the classical Cartan triality
through the morphism $\cS_4\ra\cS_3$ induced by the permutation 
of the three partitions of four objects in two pairs. In terms 
of representations, this translates into the permutation of the 
three non equivalent eight dimensional representations of $\fspin_8$:
\begin{eqnarray}
\nonumber \Delta_1 & = & (A_1\ot A_2)\op (A_3\ot A_4), \\
\nonumber \Delta_2 & = & (A_1\ot A_3)\op (A_2\ot A_4), \\
\nonumber \Delta_3 & = & (A_1\ot A_4)\op (A_2\ot A_3). 
\end{eqnarray}

\subsection*{Binary, ternary, and triality models}

Among the models of the exceptional Lie algebras that we 
will meet in the sequel, most will be derived from the 
{\it triality model} first defined, in a more general context,
by Allison in \cite{all}, and rediscovered in \cite{LMadv}. 
The idea is to associate to a (complexified)
real normed algebra $\AA=\RR,\CC,\HH,\OO$, its 
{\it triality algebra} $\ft(\AA)$ with its three natural 
modules $\AA_1,\AA_2,\AA_3$. Then for any pair $\AA,\BB$ of 
normed algebras, the direct sum
$$\fg(\AA,\BB)=\ft(\AA)\times\ft(\BB)\op (\AA_1\ot\BB_1)
\op (\AA_2\ot\BB_2)\op (\AA_3\ot\BB_3)$$
 has a natural Lie algebra structure. This leads to the 
famous Freudenthal Magic Square, whose fourth line 
$\fg(\AA,\OO)$ is the series of exceptional Lie algebras
$\ff_4, \fe_6, \fe_7, \fe_8$.  

The Lie algebras on the second and third lines of the Magic
Square are endowed each with a special module: $\fg(\AA,\CC)$ with
the cubic Jordan algebra $J_3(\AA)$, and $\fg(\AA,\HH)$ with 
the Zorn algebra $\fz_2(\AA)$. The natural inclusions 
$\fg(\AA,\CC)\subset\fg(\AA,\HH)\subset\fg(\AA,\OO)$ then 
lead to the {\it binary} and {\it ternary} models for the 
exceptional Lie algebras:
\begin{eqnarray}
\nonumber \fg(\AA,\OO) & = & \fsl_2\times\fg(\AA,\HH)
\op (\CC^2\ot \fz_2(\AA)), \\
\nonumber \fg(\AA,\OO) & = & \fsl_3\times\fg(\AA,\CC)
\op (\CC^3\ot J_3(\AA))\op (\CC^3\ot J_3(\AA))^*.
\end{eqnarray}
This extends to $\fspin_8=\ft(\OO)$, whose ternary
model is the four-ality model related to the Reye configuration. 
Note also that the triality models can be interpreted as
$\HH$-graded Lie algebras, with a similar definition to the one we 
introduced for $\OO$-gradings. 

\section{Lines on the cubic surface}

The configuration of the 27 lines on a smooth cubic surface 
in $\CC\PP^3$ have been thoroughly investigated by the classical
algebraic geometers. It has been known for a long time that 
the automorphism group of this configuration can be identified
with the Weyl group of the root system of type $E_6$, of order
$51\, 840$ \cite{cox}. Moreover, the minimal representation $J$ of the simply 
connected complex Lie group of type $E_6$ has dimension $27$. 
This is a minuscule representation, meaning that the weight
spaces are lines and that the Weyl group $W(E_6)$ acts transitively 
on the weights. In fact one can recover the lines configuration
of the cubic surface by defining two weights to be incident 
if they are not orthogonal with respect to the unique (up to 
scale) invariant scalar product. 

Conversely, one can recover the action of the Lie group $E_6$
on $J$ from the line configuration. Faulkner defines a cubic 
form on $J$ as the sums of signed monomials associated to
the tritangent planes \cite{faul}. The stabilizer of that cubic form in 
$GL(J)$ is precisely $E_6$. Note that the polarization of 
this cubic form is a symmetric bilinear map $J\times J\ra J^*$. 
Identifying appropriately 
$J$ with $J^*$ we get an algebra structure which is 
known to co\"{\i}ncide with the exceptional complex Jordan algebra $J_3(\OO)$. 

The closed $E_6$-orbit in the projectivization $\PP J_3(\OO)$
is known as the complex {\it Cayley plane} $\OO\PP^2$
and should be thought of
as the projective plane over the Cayley algebra of octonions. 
Being the orbit of any
weight space it is circumscribed to the Schoute polytope $2_{21}$, 
which appears as a discrete version of the Cayley plane. 
In particular the $10$ lines incident to a given line correspond 
to the polar quadric or $\OO$-line (whose Euler characteristic is $10$). 
The property that two general $\OO$-lines on the Cayley plane 
have a unique intersection point, thus mirrors the obvious  
fact that, two concurrent lines on the cubic surface being given, 
there exists a unique line meeting both.

\smallskip
It is well known that most of the interesting subgroups of $W(E_6)$ can 
be realized as stabilizers of some subconfigurations. 
It seems not to have been noticed before that most of them also
have natural interpretations in terms of branching. By this we
mean that we can find a subalgebra of $\fe_6$ such that the restriction
of the representation in $J$ splits in such a way that the 
relevant subconfiguration can immediately be read off. 

There is a general recipe to identify semisimple subalgebras of 
a simple complex Lie algebra, that we illustrate with the case 
of $\fe_6$ (see \cite{ov}, Chapter 6). 
One begins with the affine Dynkin diagram, which in the case
we are interested in has a remarkable threefold symmetry:

\begin{center}
\setlength{\unitlength}{3mm}
\begin{picture}(30,7)(-10,-5)

\multiput(0,0)(2,0){5}{$\circ$}
\multiput(0.5,.4)(2,0){4}{\line(1,0){1.6}} 
\put(4,-2){$\circ$}
\put(4,-4){$\circ$}
\multiput(4.3,-3.4)(0,2){2}{\line(0,1){1.5}}

\end{picture} 
\end{center}

Then we choose a set of nodes, that we mark in black. Suppressing
these nodes we get the Dynkin diagram (usually disconnected) of 
a semisimple Lie subalgebra  $\fh$ of $\fe_6$
which is uniquely defined up to conjugation. 
The Weyl group $W$ of this semisimple Lie algebra is a subgroup of $W(E_6)$,
also uniquely defined up to conjugation. We get three types of data:
\begin{enumerate}
\item {\it Combinatorial data}: $W$ can be realized as the stabilizer 
of a certain subconfiguration of the $27$ lines, encoded
in the marked Dynkin diagram;
\item {\it Representation theoretic data}: as an $\fh$-module, $J$ 
splits into a direct sum of irreducible components;
\item {\it Geometric data}: the $\fh$-components encode certain 
special subvarieties of the Cayley plane. 
\end{enumerate}

\medskip\noindent
{\it Example 1}. We mark two of the three extreme nodes. Then 
$\fh=\fspin_{10}$ and $W=W(D_5)=\ZZ_2^4\rtimes\cS_5$. 

\begin{center}
\setlength{\unitlength}{3mm}
\begin{picture}(30,7)(-10,-5)

\multiput(0,0)(2,0){5}{$\circ$}
\multiput(0,0)(8,0){2}{$\bullet$}
\multiput(0.5,.4)(2,0){4}{\line(1,0){1.6}} 
\put(4,-2){$\circ$}
\put(4,-4){$\circ$}
\multiput(4.3,-3.4)(0,2){2}{\line(0,1){1.5}}

\end{picture} 
\end{center}

The index of $W$ in $W(E_6)$ is $27$: this subgroup is just the
stabilizer of some line in the configuration. In fact $\fh$ is the 
semisimple part of the Lie algebra of 
the stabilizer of a one-dimensional weight 
space $\ell$, which defines a 
point on the Cayley plane and can be identified with one 
of the lines of the configuration. The
branching, i.e. the decomposition of $J$ as an $\fh$-module, gives 
$$J=\ell\op\Delta\op U.$$ 

The $16$-dimensional half-spin representation $\Delta$ can be 
identified with the
tangent space to the Cayley plane at $\ell$; combinatorially, 
the sixteen weight spaces generating $\Delta$ give the sixteen lines 
which do not meet $\ell$; geometrically, 
the intersection of the Cayley plane with its tangent space 
at $\ell$ is a cone over a ten dimensional spinor variety.

The $10$-dimensional natural representation $U$ encodes the normal
space to the Cayley plane at $\ell$; combinatorially, 
the ten weight spaces generating $U$ give the ten incident 
lines to $\ell$. Note that this representation is self-dual,
so its weights occur in opposite pairs corresponding to 
incident pairs of incident lines to $\ell$.
Geometrically, 
the intersection of the Cayley plane with $\PP U$ is 
the polar eight-dimensional quadric, a copy of the 
projective line $\OO\PP^1$ over the Cayley algebra. 

\medskip\noindent
{\it Example 2}. We mark the three extreme nodes. In this case
$\fh=\fspin_8$ and $W=W(D_8)=\ZZ_2^3\rtimes\cS_4$. 

\begin{center}
\setlength{\unitlength}{3mm}
\begin{picture}(30,7)(-10,-5)

\multiput(0,0)(2,0){5}{$\circ$}
\multiput(0,0)(8,0){2}{$\bullet$}
\multiput(0.5,.4)(2,0){4}{\line(1,0){1.6}} 
\put(4,-2){$\circ$}
\put(4,-4){$\bullet$}
\multiput(4.3,-3.4)(0,2){2}{\line(0,1){1.5}}

\end{picture} 
\end{center}

By restricting the previous case we get the branching 
$$J=\ell_1\op\ell_2\op\ell_3\op\Delta_1\op\Delta_2\op\Delta_3,$$
where $\Delta_1,\Delta_2,\Delta_3$ are the three eight dimensional
representations of $Spin_8$, which we deliberately avoid to 
distinguish since they are exchanged by Cartan's triality. 
The three lines $\ell_1,\ell_2,\ell_3$ are pairwise incident, 
hence they are the three intersection lines of the cubic
surface with a {\it tritangent plane}. Note that the 
index of $W$ in $W(E_6)$ is $270=6\times 45$. Since we have a sixfold
ambiguity on the order of the three lines, we recover the 
classical fact that the cubic surface has exactly $45$
tritangent planes.  

Once we have fixed these three lines, each of them has eight more
incident lines coming into four pairs, and corresponding to the 
pairs of opposite weights of one of the eight-dimensional 
representations of $\fspin_8$. 
Note that this exhausts the $27$ lines. 

\begin{center}
\setlength{\unitlength}{3mm}
\begin{picture}(20,14)(-2,-.5)
\put(0,0){ \resizebox{!}{1.5in}{\includegraphics{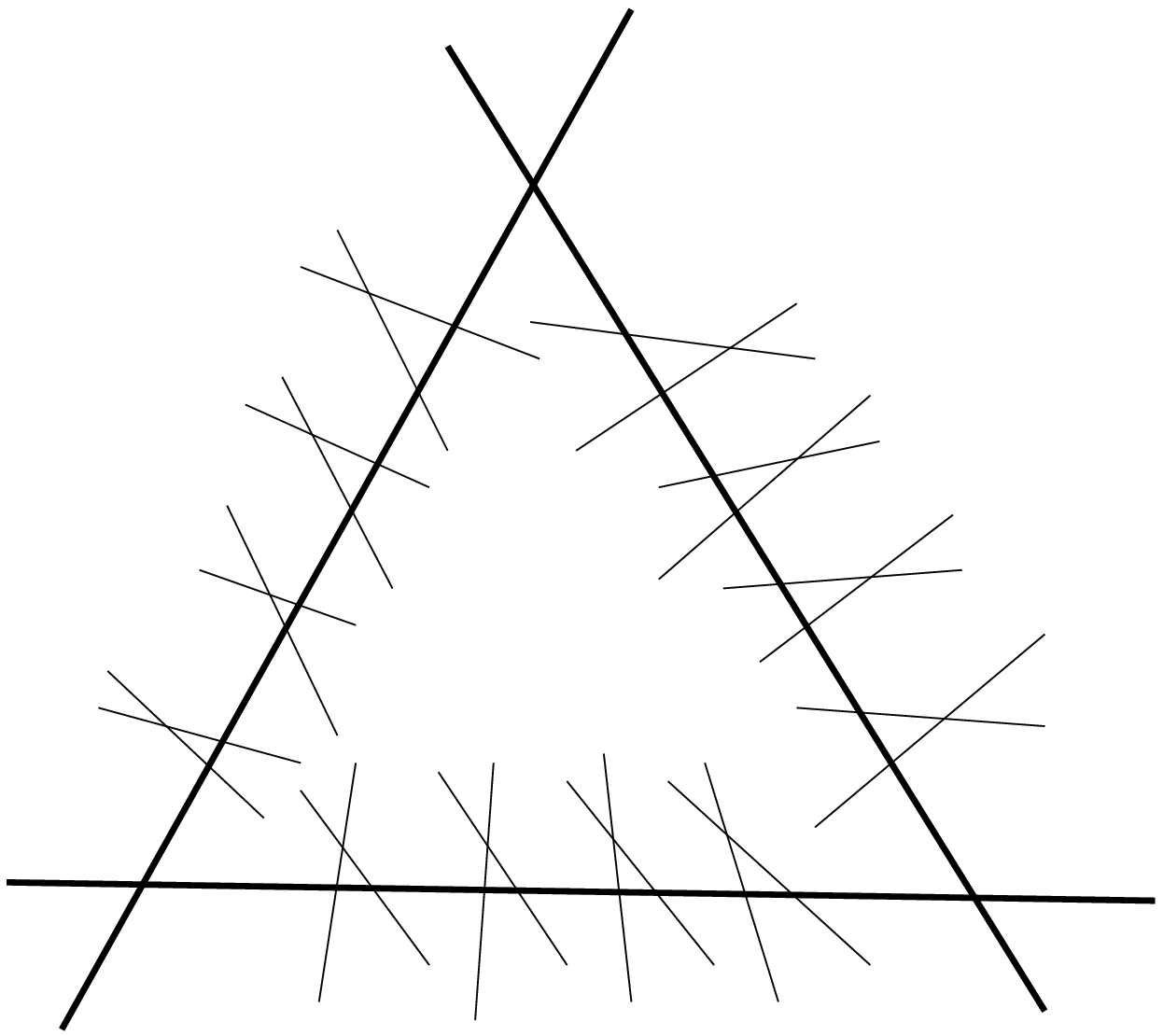}}}
\put(9,12){$\ell_1$}
\put(15,1){$\ell_2$}
\put(4.5,12){$\ell_3$}
\end{picture} 
\end{center}

\centerline{{\it Figure 3. Triality from the $27$ lines}}

\medskip
The sum of the weights of $\ell_1,\ell_2,\ell_3$ is zero,
and this characterizes triples of lines on a tritangent plane. 

Geometrically, we have three eight-dimensional quadrics on the
Cayley plane, any two of them meeting exactly in one point. 
In terms of plane projective geometry, these three quadrics 
are projective lines which are the sides of a self-polar triangle.

\medskip\noindent
{\it Example 3}. Now we mark a unique node, which is neither extremal
nor central. Then
$\fh=\fsl_2\times \fsl_6$ and $W=\cS_2\times\cS_5$. This leads to the 
binary model of $\fe_6$.

\begin{center}
\setlength{\unitlength}{3mm}
\begin{picture}(30,7)(-10,-5)

\multiput(0,0)(2,0){5}{$\circ$}
\multiput(0.5,.4)(2,0){4}{\line(1,0){1.6}} 
\put(4,-2){$\bullet$}
\put(4,-4){$\circ$}
\multiput(4.3,-3.4)(0,2){2}{\line(0,1){1.5}}

\end{picture} 
\end{center}

The index of $W$ in $W(E_6)$ is $36$. The branching gives
$$J=U\otimes A\op\Lambda^4U,$$ 
where $U$ denotes the six-dimensional natural representation 
of $\fsl_6$, and $A$ the natural representation of $\fsl_2$. 
The twelve weights of $U\otimes A$ split into six pairs
$(\ell_i,\ell'_i)$ where $\ell_1,\ldots ,\ell_6$ have 
the same component over $A$. Of course 
$$\begin{pmatrix} \ell_1 & \ell_2 & \ell_3 & \ell_4 & \ell_5 & \ell_6 \\
\ell'_1 & \ell'_2 & \ell'_3 & \ell'_4 & \ell'_5 & \ell'_6 
\end{pmatrix}$$
is a double-six, and there are exactly $36$ double-sixes on the
cubic surface. 

Note that the node that we have marked is the node of the Dynkin 
diagram of $E_6$ that defines the adjoint representation. This
explains the correspondence between double-sixes and pairs of 
opposite roots of $E_6$ (\cite{dolg}, 10.1.5).  

Geometrically, such a pair of roots defines a point of the adjoint
variety, the projectivization of the minimal nilpotent orbit in
the Lie algebra $\fe_6$. The image of its action on $J$ has minimal
rank, namely $6$ \cite{IM}, and is a maximal linear space in the Cayley plane, 
with a weight basis given by one half of the double-six. The other
half can be recovered from the similar action on $J^*$, whose image 
is again a maximal linear space in (a dual copy of) the Cayley plane. 

Explicitely, let $P(J)$ denote 
the set of weights of $J$ (the $W(E_6)$-orbit of the fundamental 
weight $\omega_1$, in the notation of \cite{bou}). 
Then the double-six $D_{\a}$ 
associated to a root $\a$ (up to sign), considered as a set of
such weights, is
$$\begin{array}{rcl}
D_{\a}& = & \{\g\in P(J),\; \g+\a\; or\; \g-\a\in P(J)\} \\
 & =& \{\g\in P(J),\; (\g,\a)\ne 0\}.
\end{array}$$ 

Now given two double-sixes $D_{\a}$ and $D_{\b}$, they can have 
only two different relative positions, following that $(\a,\b)=0$
or not.  In the first case they are said to be {\it azygetic}, 
and $\#(D_{\a}\cap D_{\b})=6$. Then $\a+\b$ or $\a-\b$ is again a root 
and defines a third double-six $D_{\a\pm\b}$, which  is azygetic
to both $D_{\a}$ and $D_{\b}$. There exist $120$ such {\it azygetic 
triads of double-sixes}, corresponding to the $120$ subsystems of type
$A_2$ of the root system $E_6$. 

In the latter case the double-sixes are {\it syzygetic},
and $\#(D_{\a}\cap D_{\b})=4$. As indicated by the 
marked Dynkin diagram

\begin{center}
\setlength{\unitlength}{3mm}
\begin{picture}(30,7)(-10,-5)

\multiput(0,0)(2,0){5}{$\circ$}
\multiput(0.5,.4)(2,0){4}{\line(1,0){1.6}} 
\put(4,-2){$\bullet$}\put(2,0){$\bullet$}\put(6,0){$\bullet$}
\put(4,-4){$\circ$}
\multiput(4.3,-3.4)(0,2){2}{\line(0,1){1.5}}

\end{picture} 
\end{center}

\noindent we can complete such a pair of double-sixes into a 
{\it syzygetic tetrad} of double-sixes. If we let $A_1, A_2, 
A_3, A_4$ be the two-dimensional natural representations of the 
four copies of $\fsl_2$ corresponding to the white nodes of the 
diagram, the minimal representation of $\fe_6$ branches to
$$J=\CC^3\op \bigoplus_{i<j}A_i\ot A_j,$$
and the associated tetrad of double-six is given by the weights 
of the four submodules 
$$D_i=\bigoplus_{j\ne i} A_i\ot A_j.$$
The number of syzygetic tetrads is the number of root subsystems
of type $A_1^4$ in $E_6$, that is $135$. 

\medskip\noindent
{\it Example 4}. Now we mark the central node. Then
$\fh=\fsl_3\times \fsl_3\times \fsl_3$ and $W=\cS_3\times\cS_3\times\cS_3$.
This leads to the ternary model of $\fe_6$. 

\begin{center}
\setlength{\unitlength}{3mm}
\begin{picture}(30,7)(-10,-5)

\multiput(0,0)(2,0){5}{$\circ$}
\multiput(0.5,.4)(2,0){4}{\line(1,0){1.6}} 
\put(4,0){$\bullet$}
\put(4,-2){$\circ$}
\put(4,-4){$\circ$}
\multiput(4.3,-3.4)(0,2){2}{\line(0,1){1.5}}

\end{picture} 
\end{center}

The index of $W$ in $W(E_6)$ is $240$. The branching gives
$$J=(A^*\otimes B)\op (B^*\otimes C)\op (C^*\otimes A),$$ 
where $A,B,C$ denote the natural representations
of the three copies of $\fsl_3$ in $\fh$. The $27$ weights
are thus split into three bunches of nine. 

Consider for example
the nine weights $\epsilon'_i-\epsilon_j$ of $A^*\otimes B$, where the
$\epsilon_j$ are the weights of $A$ and the $\epsilon'_i$ those of
$B$ (note that both sets sum to zero). 
Display these weights on a $3\times 3$ square by putting 
$\epsilon'_{i+j-1}
-\epsilon_{i+2j-2}$ in the box $(i,j)$, the indices being 
considered modulo three. 
Then three weights on the same line or column have they sum 
equal to zero, and thus define a tritangent plane. We have 
obtained what is called a {\it Steiner set} -- nine lines 
obtained as the intersections of two trihedra. 

Moreover, we have split the $27$ lines into three such sets forming 
a so-called {\it Steiner triple system}. Steiner sets are in
correspondence with root subsystems of type $A_2$ of the root system
$E_6$. The orthogonal to such a subsystem is the orthogonal 
product of two other $A_2$-subsystems. This is why, a Steiner set 
being given, there is a unique way to complete it into a triple
system (\cite{dolg}, 10.1.7). Of course this is another manifestation 
of triality (see \cite{edge3})!

Note also that the invariant cubic form on $J$ can be characterized, up
to scalar, as the unique $W(E_6)$-invariant cubic form whose restriction 
to each factor of type $A^*\ot B$ is proportional to the determinant. 

Geometrically, three Steiner sets correspond to
 three $\PP^8$'s in $\PP J$, each  cutting
the Cayley plane along a copy of the Segre variety $\PP^2\times\PP^2$. 

\medskip\noindent {\it Remark.}
Define two Steiner sets to be incident if they have no common 
tritangent plane. Each Steiner set is incident to exactly 
$56$ other sets, including the two special ones which complete
it into a triple system. If we choose one of these, exactly 
$28$ Steiner sets are incident to both, including the remaining
set in the triple system. Contrary to what we could be tempted to
believe, the configuration of the remaining $27$ Steiner sets is not 
combinatorially equivalent to that of the $27$ lines. Indeed, one
can check that each of the $27$ Steiner sets is incident to only
eight of the other ones. Nevertheless, the Steiner sets define
an interesting regular graph, which has the same number of edges 
and vertices than the graph defined by the diameters of the polytope
$4_{21}$, although it is not combinatorially equivalent. 

Coxeter already noticed in \cite{coxeter}
that the $40$ triple systems can be interpreted as 
hexagons on the polytope $4_{21}$

\section{Bitangents to the plane quartic curve}
 
The $28$ bitangents to a smooth plane quartic curve give rise to 
$56$ lines on the Del Pezzo surface of degree two defined
as the double cover of the projective plane, branched over 
the quartic \cite{dem}. This line configuration has for automorphism 
group the index two normal subgroup $W(E_7)^+$ of the 
Weyl group of $E_7$, which has order $2\, 903\, 040=2^{10}\times 
3^{4}\times 5\times 7$. 

The Lie group of type $E_7$ has a minimal representation 
$V$ of dimension $56$, which is again minuscule. The invariant 
forms are a symplectic form -- so that the $56$ weights split
into $28$ pairs of opposite weights --  and a quartic form
which cannot be deduced, contrary to the case of the $27$
lines, solely from the configuration. The weights of this 
representation form the {\it Gosset polytope}
$3_{21}$. This polytope appears as a discrete version of the minimal 
$E_7$-orbit in $\PP V$, which we call the {\it Freudenthal variety}. 

\smallskip
The affine 
Dynkin diagram of type $E_7$ has a two fold symmetry:

\begin{center}
\setlength{\unitlength}{3mm}
\begin{picture}(30,5)(-10,-3)

\multiput(0,0)(2,0){7}{$\circ$}
\multiput(0.5,.4)(2,0){6}{\line(1,0){1.6}} 
\put(6,-2){$\circ$}
\put(6.3,-1.4){\line(0,1){1.5}}

\end{picture} 
\end{center}

Again we deduce classical configurations of bitangents from
markings of this diagram.

\medskip\noindent
{\it Example 1}. We mark the two opposite extreme nodes.
Then $\fh=\fe_6$, and $W=W(E_6)$ has index $28$ in $W(E_7)$. 

\begin{center}
\setlength{\unitlength}{3mm}
\begin{picture}(30,5)(-10,-3)

\multiput(2,0)(2,0){5}{$\circ$}
\multiput(0,0)(12,0){2}{$\bullet$}
\multiput(0.5,.4)(2,0){6}{\line(1,0){1.6}} 
\put(6,-2){$\circ$}
\put(6.3,-1.4){\line(0,1){1.5}}

\end{picture} 
\end{center}

Indeed
it is well-known that the stabilizer of any bitangent is a copy 
of the automorphism group of the $27$ lines. Since the action 
of the latter is irreducible, the branching has to give an irreducible 
supplement, up to the sign of the weights. And indeed we have the 
decomposition into $\fh$-modules
$$V=\CC\op J\op J^*\op\CC.$$ 
Geometrically, the factor $J=J_3(\OO)$ appears as the tangent space to the 
Freudenthal variety, and $J^*$ is the first normal space. 
Note that the intersection of the Freudenthal variety with its 
tangent space is a cone over the Cayley plane, whose discrete 
skeleton is precisely given by the $27$ lines. 

\medskip\noindent
{\it Example 2}. Mark only one node next to one of the
opposite extremal nodes. Then $\fh=\fsl_2\times \fspin_{12}$, and 
$W$ has index $63$ in $W(E_7)$. This leads to the binary 
model of $\fe_7$.

\begin{center}
\setlength{\unitlength}{3mm}
\begin{picture}(30,5)(-10,-3)

\multiput(0,0)(2,0){7}{$\circ$}
\put(2,0){$\bullet$}
\multiput(0.5,.4)(2,0){6}{\line(1,0){1.6}} 
\put(6,-2){$\circ$}
\put(6.3,-1.4){\line(0,1){1.5}}

\end{picture} 
\end{center}

Let $A, W$ denote the natural representations 
of $\fsl_2$ and $\fspin_{12}$. Here the branching gives 
the very simple decomposition
$$V=A\ot W\op \Delta,$$
where $\Delta$ is one of the half-spin representations, of dimension
$32$. The factor $A\ot W$ corresponds to a set of twelve bitangents. 
Since $W$ has an invariant quadratic form, its weights come into 
pairs of opposite weights. We thus get six pairs of bitangents 
forming a {\it Steiner complex} (\cite{dolg}, 6.1.2). 

Since the node that we have marked defines the adjoint representation 
of $\fe_7$ on the Dynkin diagram of type $E_7$, the $63$ Steiner 
complexes are in natural correspondence with the $63$ pairs of
opposite roots in the root system $E_7$. From this 
perspective they play the same role as the double-sixes of lines 
on the cubic surface.

Geometrically, the twelve weight spaces of a Steiner complex
generate a $\PP^{11}$ in $\PP W$ whose intersection with the
Freudenthal variety is a ten dimensional quadric. If we denote by 
$P(V)$ the set of weight of $V$ (the $W(E_7)$-orbit of the 
fundamental weight $\omega_1$, in the notation of \cite{bou}),
the Steiner complex associated to the root $\a$ (up to sign) is 
$$S_{\a}=\{\g\in P(V),\; \g+\a\; {\mathrm or}\; \g-\a\in P(V)\}.$$ 
Now if $\b\ne\pm\a$ is another root, we can have $(\a,\b)=0$ or not,
with respectively $30$ and $32$ possibilities for $\b$, up to sign.
 
In the first case, the corresponding root spaces generate in $\fe_7$
a copy of $\fsl_2\times\fsl_2$, and the Steiner complexes are 
{\it syzygetic}, which means that $\#(S_{\a}\cap S_{\b})=4$.
The roots which are orthogonal both to $\a$ and $\b$ form a 
reducible root system of type $A_1\times D_4$. In particular, 
there is a third Steiner complex $S_{\g}$ canonically associated 
to the pair $(S_{\a},S_{\b})$, and syzygetic to both of them. 
The characteristic property of the triple  $(S_{\a},S_{\b},S_{\g})$
is that its union is the whole set of bitangents (see Example 6 below). 
The number of such {\it syzygetic triads of Steiner complexes} is 
$63\times 30/6=315$. 

In the second case, the root spaces generate 
a copy of $\fsl_3$, and the Steiner complexes are 
{\it azygetic}: $\#(S_{\a}\cap S_{\b})=6$. 
Then $\a+\b$ (or $\a-\b$) is again a root, and $S_{\a+\b}$ 
is azygetic both to $S_{\a}$ and $S_{\b}$: we obtain an 
{\it azygetic triad of Steiner complexes}. The number of such
triads is $63\times 32/6=336$.

\medskip\noindent
{\it Example 3}. We mark the lowest node.
Then $\fh=\fsl_8$, and $W=\cS_8$ has index $36$ in $W(E_7)$. 

\begin{center}
\setlength{\unitlength}{3mm}
\begin{picture}(30,5)(-10,-3)

\multiput(0,0)(2,0){7}{$\circ$}
\multiput(0.5,.4)(2,0){6}{\line(1,0){1.6}} 
\put(6,-2){$\bullet$}
\put(6.3,-1.4){\line(0,1){1.5}}

\end{picture} 
\end{center}

If we denote by $U$ the natural representation of $\fsl_8$, the branching
gives 
\begin{eqnarray}
\nonumber V &= &\Lambda^2U\op\Lambda^2U^*, \\
\nonumber \fe_7 &= &\fsl_8\op\Lambda^4U.
\end{eqnarray}
Once we fix a basis $u_1,\ldots ,u_8$ of $U$, we can therefore
identify each bitangent with a pair $(ij)$, with $1\le i<j\le 8$. 
Such a notation seems to have been first introduced by Hesse.
Moreover, the Weyl group of $E_7$ is generated by the symmetric
group $\cS_8$ and the symmetries associated to the roots coming from 
the factor $\Lambda^4U$. These symmetries are indexed by partitions
$(pqrs|xyzt)$ of $[1,8]$ into disjoints fourtuples. The induced action 
on the bitangents is given by 
$$s_{(pqrs|xyzt)}(ij)=\left\{  \begin{array}{ll} 
(pqrs/ij) & {\mathrm if}\quad \{ij\}\subset\{pqrs\}, \\
(xyzt/ij) & {\mathrm if}\quad \{ij\}\subset\{xyzt\}, \\
(ij) & {\mathrm otherwise}.
\end{array}    \right.$$
This is classicaly called a {\it bifid transformation}.

\medskip\noindent
{\it Example 4}. We mark two extreme but not opposite nodes.
Then $\fh=\fsl_7$, and $W=\cS_7$ has index $288$ in $W(E_7)$. 

\begin{center}
\setlength{\unitlength}{3mm}
\begin{picture}(30,5)(-10,-3)

\multiput(2,0)(2,0){6}{$\circ$}
\put(0,0){$\bullet$}
\multiput(0.5,.4)(2,0){6}{\line(1,0){1.6}} 
\put(6,-2){$\bullet$}
\put(6.3,-1.4){\line(0,1){1.5}}

\end{picture} 
\end{center}

Obviously this example comes from the previous one: we have 
just passed from $\fsl_8$ to $\fsl_7$. If we denote by $U$ the 
natural representation of $\fsl_7$, the branching gives 
$$V=U\op\Lambda^2U\op\Lambda^2U^*\op U^*.$$
We have thus distinguished a set of seven bitangents forming
what is called an {\it Aronhold set} (\cite{dolg}, 6.1.3). 
Geometrically, the seven
weight spaces of an Aronhold set generate a $\PP^6$ which is
a maximal linear space inside the Freudenthal variety. 

Note that the Aronhold sets generate the remaining basic
representation $R$ of $\fe_7$, the one corresponding to the lowest 
extremal node of the Dynkin diagram. By this we mean that among
the lines of $\Lambda^7V$ generated by the weight vectors forming
Aronhold sets, one has a dominant weight and generates a copy of 
$R$. We know that any irreducible representation of $\fe_7$ can 
then be constructed from $V$, $R$ and the adjoint representation 
by natural tensorial operations.

\medskip\noindent
{\it Example 5}. We mark the central node. 
Then $\fh=\fsl_4\times \fsl_4\times \fsl_2$, and $W=\cS_4\times 
\cS_4\times \cS_2$ has index $1260$ in $W(E_7)$. 

\begin{center}
\setlength{\unitlength}{3mm}
\begin{picture}(30,5)(-10,-3)

\multiput(0,0)(2,0){7}{$\circ$}
\put(6,0){$\bullet$}
\multiput(0.5,.4)(2,0){6}{\line(1,0){1.6}} 
\put(6,-2){$\circ$}
\put(6.3,-1.4){\line(0,1){1.5}}

\end{picture} 
\end{center}

Denote by $C$ the two-dimensional natural representation,
by $A,B$ the two four dimensional ones. Then the branching gives
$$V=(\Lambda^2A\op\Lambda^2B)\otimes C \op (A\ot B\op A^*\ot B^*).$$
Note that $\Lambda^2A$ and $\Lambda^2B$ are self-dual, as well as
$C$, so we have distinguished two sets of $6$ bitangents forming a
Steiner complex. 

The remaining sixteen bitangents are indexed by the weights 
of $A\ot B$. Recall that four bitangents form a {\it syzygetic 
tetrad} when their eight tangency points are the eight  
intersection points of the plane quartic with some conic
(\cite{dolg}, 6.1.1). In
terms of weights, this means that the four bitangents can 
be represented by weights summing to zero. Here we observe
a phenomenon very similar to the property of a Steiner set of 
lines on the cubic surface: our sixteen bitangents can be 
split into four syzygetic tetrads in essentially twelve
different ways. Indeed our tetrads must be of the form 
$$T_{\s}=\{\epsilon_i+\epsilon_{\s(i)},\quad 1\le i\le 4\}$$
for some permutation $\s$, and we have to find four permutations
$\s_1$, $\s_2$, $\s_3$, $\s_4$ such that $\s_j(i)\ne \s_k(i)$
for each $i$ and each $j\ne k$. The twelve possibilities  are given
by the four-tuples of permutations of the form 
$(pqrs)$, $(qpsr)$, $(rspq)$, $(srqp)$ or 
$(pqrs)$, $(qpsr)$, $(rsqp)$, $(srpq)$.  

\medskip\noindent
{\it Example 6}. Now we mark a non extremal node next to the central
one.  Then $\fh=\fsl_3\times \fsl_6$, and $W=\cS_3\times 
\cS_6$ has index $336$ in $W(E_7)$. This case leads to the ternary 
model of $\fe_7$.

\begin{center}
\setlength{\unitlength}{3mm}
\begin{picture}(30,5)(-10,-3)

\multiput(0,0)(2,0){7}{$\circ$}
\put(4,0){$\bullet$}
\multiput(0.5,.4)(2,0){6}{\line(1,0){1.6}} 
\put(6,-2){$\circ$}
\put(6.3,-1.4){\line(0,1){1.5}}

\end{picture} 
\end{center}

Let $A, B$ denote the natural representations of $\fsl_3$ and $\fsl_6$. 
The branching gives the decomposition 
$$V=A\ot B \op \Lambda^3B\op A^*\ot B^*,$$
where the middle factor is self-dual. The factor $A\ot B$ splits, 
following the $A$-component, into three sixers of bitangents. 
Aggregating them  two by two we get
an azygetic triad of Steiner sets. 

\medskip\noindent
{\it Example 7}. We mark the two nodes next to the two 
opposite extremal nodes. Then $\fh=\fsl_2\times \fsl_2\times\fspin_8$, and 
$W$ has index $2\times 315$ in $W(E_7)$. 

\begin{center}
\setlength{\unitlength}{3mm}
\begin{picture}(30,5)(-10,-3)

\multiput(0,0)(2,0){7}{$\circ$}
\multiput(2,0)(8,0){2}{$\bullet$}
\multiput(0.5,.4)(2,0){6}{\line(1,0){1.6}} 
\put(6,-2){$\circ$}
\put(6.3,-1.4){\line(0,1){1.5}}

\end{picture} 
\end{center}

Let $A, B$ denote the natural representations of the two copies 
of $\fsl_2$. Up to triality we may suppose that the natural (non
spinorial) eight-dimensional representation $W$ of $\fspin_8$ is given
by the lowest node.
Then the branching gives the decomposition 
$$V=A\ot\Delta_+\op B\ot\Delta_-\op 2W\op 2A\ot B,$$
where the last two factors are two copies of the same module,
put on duality by the symplectic form.  

The symmetry of the picture suggests that we should introduce a 
supplementary copy of $\fsl_2$, hence three copies with natural
representations $A_1, A_2, A_3$ such that 
$$V=A_1\ot\Delta_1\op A_2\ot\Delta_2\op A_3\ot\Delta_3\op
A_1\ot A_2\ot A_3.$$
Indeed this is precisely what the trialitarian description of 
$\fe_7$ tells us (see \cite{LMadv}, Theorem 4.1). 
We get a partition of the $28$ bitangents 
into three groups of eight and a
group of four. Adding the latter to the three former we get 
a {\it syzygetic triad of Steiner complexes}. 

We have exactly $315$ such triads, and this is also 
the number of syzygetic tetrads (\cite{dolg}, Corollary 6.1.4). 
Indeed the weights of 
$A_1\ot A_2\ot A_3$ define such a tetrad.

We recapitulate:

\begin{prop}
There are natural correspondences between:
\begin{enumerate}
\item Steiner complexes of bitangents and root subsystems of type 
$A_1$; 
\item azygetic triples of Steiner complexes and subsystems of type 
$A_2$;
\item syzygetic triples of Steiner complexes and root subsystems of type 
$D_4$ of the root system $E_7$.
\end{enumerate}
\end{prop}

\subsection*{Symplectic geometry}
We have already mentionned that the Weyl group of $E_7$ 
can (almost) be identified with a classical group over a finite field 
(see \cite{bou}, Exercice 3 of section 4, p. 229), 
namely $$W(E_7)^+\simeq Sp(6,\FF_2).$$
This means that the incidence geometry of the $28$ bitangents
should be interpreted as a symplectic six-dimensional geometry 
over the field with two elements. Such symplectic geometries 
appear on the third line of Freudenthal's magic square, and the 
$28$ bitangents give a finite model of these geometries (see e.g.
\cite{LMpop}). 

Recall that a symplectic five-dimensional projective geometry 
has three types of elements: points, isotropic lines and isotropic 
planes. In the $E_7$ geometry, points and isotropic lines are points 
and lines in the Freudenthal variety, while isotropic planes are
in correspondence with maximal, ten dimensional quadrics on 
the Freudenthal variety. 

In our finite geometry,
we have seen that the points correspond to the $63$ 
Steiner complexes. The $315$ isotropic lines are in 
correspondence 
with the syzygetic triads of Steiner complexes, where
we split the $28$ bitangents into three sets of twelve,
with four bitangents common to the three. It is obvious 
from that point of view that two syzygetic complexes can be 
uniquely completed into a syzygetic triad: indeed, they
are syzygetic if they generate an isotropic line, and there 
is a unique other point on that line.  

\smallskip
What are the $135$ planes? Classically, they are called G\"opel
spaces (see \cite{coble}, Chapter II, 22, and \cite{dolgort}, Chapter
IX)  and play an important role in the study of the Schottky 
problem. But let us rather skip to our representation theoretic 
point of view.
 
A projective plane over $\FF_2$ 
is a Fano plane. It has $7$ points 
and $7$ lines. So to get an isotropic plane 
we should be able to partition the $28$ 
bitangents into $7$ quadruples, in such a way that the 
complement of each of them can be split into three octuples
defining syzygetic Steiner complexes. This looks like a
combinatorial challenge but representation theory tell us 
what to do! We have already used the fact that $\fe_7$ has 
a maximal semisimple Lie algebra isomorphic to $\fso_8\times\fsl_2^3$. 
Thanks to the four-ality model of $\fso_8$ we can take 
four two-dimensional spaces $A_4, A_5, A_6, A_7$ and decompose 
$$\fso_8=\fsl(A_4)\times\fsl(A_5)\times\fsl(A_6)\times\fsl(A_7)
\op (A_4\ot A_5\ot A_6\ot A_7).$$
Then the three eight-dimensional representations decompose as 
\begin{eqnarray}
\nonumber \Delta_1 & = & (A_4\ot A_5)\op (A_6\ot A_7), \\
\nonumber \Delta_2 & = & (A_4\ot A_6)\op (A_5\ot A_7), \\
\nonumber \Delta_3 & = & (A_4\ot A_7)\op (A_5\ot A_6). 
\end{eqnarray}
Then we plug that in the decomposition of the $56$-dimensional
representation of $E_7$ given in Example 7 above. 
We obtain, if we denote $A_{ijk}=A_i\ot
A_j\ot A_k$:
$$V=A_{123}\op A_{145}\op A_{167}\op A_{246}\op 
A_{257}\op A_{347}\op A_{356}.$$
The seven Steiner complexes that we are looking for are the 
sets of weights of the submodules 
$$S_i=\bigoplus_{ijk}A_{ijk},$$
and the seven syzygetic triads they form are given by the
weights of the three submodules $S_i, S_j, S_k$ for $(ijk)$
one of the seven triples in the decomposition of $V$. 

Note that these seven triples of indices have the crucial
property that any pair of integers between
one and seven, appear in one and only one of them. Otherwise
said, they form a {\it Steiner triple system} $S(2,3,7)$ (see e.g. 
\cite{conw}). Up to isomorphism there is only one such system, 
given by the lines of the Fano plane, as one can see
on the next picture:

\begin{center}
\setlength{\unitlength}{4mm}
\begin{picture}(20,8)(-10,-1)
\put(-3.1,0){\line(1,0){6.2}} 
\put(-3.1,0){\line(3,5){3.1}} 
\put(3.1,0){\line(-3,5){3.1}} 
\put(-3,0){\line(5,3){4.44}} 
\put(3,0){\line(-5,3){4.44}} 
\put(0,0){\line(0,1){5.2}} 
\put(-2.3,2.8){2}\put(1.8,2.8){6}\put(-.3,-1){4}\put(-3.8,-1){1}
\put(3.1,-1){5}\put(-.3,5.5){3}\put(.2,2.4){7}
\put(0,1.7){\circle{4.2}}
\end{picture} 
\end{center}

In particular, by projective duality these lines can be 
represented by points of another Fano plane. 
We deduce that the stabilizer of our configuration in $W(E_7)$
is the product of $7$ copies of $A_2$ by the automorphism group
of the Fano plane, which is nothing else than the Klein group
$PSL(3,\FF_2)\simeq PSL_2(\FF_7)$, with $168$ elements. The 
index of this stabilizer is $135$, as expected. We call the 
corresponding configurations {\it Fano heptads of  Steiner
complexes}. Let us recapitulate the correspondence:

$$\begin{array}{ccc}
 {\rm Symplectic\;geometry} & \qquad {\rm Number} \qquad & {\rm Bitangents} \\
 & & \\
Points & 63 & Steiner \;complexes \\
Lines & 315 & syzygetic \;triads \\
Planes & 135 & Fano\; heptads

\end{array}$$

\medskip
\subsection*{The quartic form}
The $E_7$-module $V$ has two basic invariants tensors, such that 
$E_7$ can be described as the group of linear transformations of 
$V$ preserving these tensors: a symplectic form and a quartic form. 
The existence of such invariant forms 
 is clear on our previous decomposition of $V$. Indeed, each 
factor $A_{ijk}$ has a symplectic form induced by the choice of 
two-forms on each factor $A_l$. Moreover, $A_{ijk}$ has an invariant
quartic form given by the {\it Cayley hyperdeterminant}, which is an 
equation of the dual variety of the Segre product 
$\PP A_i\times \PP A_j\times \PP A_k\subset \PP A_{ijk}$ \cite{GKW}. 
 
\begin{prop}
The invariant quartic form on $V$ is the unique $W(E_7)$-invariant
quartic form whose restriction to each factor $A_{ijk}$ is 
proportional to the Cayley hyperdeterminant.
\end{prop}

\proof Since their is a unique invariant quartic form on a factor 
$A_{ijk}$, up to scalar, it certainly co\"{\i}ncides with the restriction 
of the invariant quartic form on $V$. Conversely, we know that up to 
the action of the Weyl group, the monomials in the invariant quartic 
form on $V$ are of three types (see \cite{lurie}): products of two, 
equal or distinct, 
products of two variables associated to opposite weights; other 
products of four variables associated to fourtuples of weights of
sum zero (thus defining syzygetic tetrads of bitangents). 
Then we must give suitable relative coefficients to the sums of 
monomials of each type. This is fixed by restriction to 
a single factor $A_{ijk}$ since the three types of monomials appear 
in the hyperdeterminant (see \cite{GKW}). \qed

\subsection*{Reconstructing $\fe_7$.}
From the trialitarian construction of $\fe_7$ and the four-ality 
for $\fso_8$ we deduce the model:

$$\fe_7=\times_{i=1}^7\fsl(A_i)\;\op\bigoplus_{(ijkl)\in I}
A_i\ot A_j\ot A_k\ot A_l,$$
where $I$ is the following set of $7$ quadruples:
$$
1247  \qquad 1256  \qquad 1346  \qquad 1357  
\qquad 2345  \qquad 2367  \qquad 4567$$
These quadruples are in natural correspondence with lines: 
simply associate to a line the four points of its complement. 
Moreover, the action on $V$ can be recovered geometrically:
each quadruple $(ijkl)$ in $I$ can be seen as a complete quadrangle in 
the Fano plane, with three pairs of opposite sides which are sent
one to the other by the $\fe_7$-action restricted to 
$A_i\ot A_j\ot A_k\ot A_l$.

\smallskip
Let us rather try to reconstruct the Lie bracket. Consider
two factors $A_i\ot A_j\ot A_k\ot A_l$ and $A_i\ot A_j\ot A_m\ot A_n$: 
we have two indices $i,j$ in common, and the third point of the line 
generated by $\a=(ijkl)$ and $\b=(ijmn)$ is $\a+\b=(klmn)$. 
The restriction of the Lie bracket defines a map
$$A_i\ot A_j\ot A_k\ot A_l\times A_i\ot A_j\ot A_m\ot A_n
 \longrightarrow  A_k\ot A_l\ot A_m\ot A_n $$
such that for some non zero constant $\theta_{\a,\b}$,
$$
 [x_i\ot x_j\ot x_k\ot x_l, y_i\ot y_j\ot y_m\ot y_n] = 
\theta_{\a,\b}\omega (x_i,y_i)\omega (x_j,y_j)
x_k\ot x_l\ot y_m\ot y_n.
$$
Indeed this is the unique equivariant map up to scalar, 
and it must be non zero because of the properties of the 
Lie bracket in a semisimple Lie algebra. 
The skew symmetry of the Lie bracket then implies that 
$$\theta_{\b,\a}=-\theta_{\a,\b}.$$
The Jacobi identity can be expressed in the following way: for
each triangle $(\a,\b,\g)$ in the Fano plane, we have the relation
$$\theta_{\a,\b}\theta_{\a+\b,\g}=\theta_{\b,\g}\theta_{\b+\g,\a}
=\theta_{\g,\a}\theta_{\g+\a,\b}.$$
The Fano plane has $28$ triangles, hence $56$ quadratic relations.

\begin{lemm}
Let $\theta_{\a,\b}=\pm 1$ according to the following rule:
the multiplication table of the canonical basis $e_1,\ldots ,e_7$ 
of the imaginary octonions is given by 
$$e_{\a}e_{\b}=\theta_{\a,\b}e_{\a+\b}\qquad for \quad \a\ne\b .$$
Then the relations above are satisfied.
\end{lemm}

\proof Denote by $e_0, e_{\a}$, where $\a=1,\ldots ,7$,
 the canonical basis of the octonions,
 (see the Introduction). Our claim amounts
to the identity 
$$(e_{\a}e_{\b})e_{\g}=(e_{\b}e_{\g})e_{\a}$$
when $\a,\b,\g$ are distinct and not aligned. To prove this we need
to remember that the Cayley algebra, although not associative, is 
{\it alternative} \cite{LMpop}. This means that the associator 
$A(x,y,z)=(xy)z-x(yz)$ is an alternating function of the arguments. 
Using the fact that $e_{\a}e_{\b}=-e_{\b}e_{\a}$ when $\a,\b$ are
distinct, we deduce that
$$\begin{array}{l}
(e_{\b}e_{\g})e_{\a}-e_{\b}(e_{\g}e_{\a})=A(e_{\b},e_{\g},e_{\a}) \\
\hspace*{3.1cm}
=-A(e_{\b},e_{\a},e_{\g})=(e_{\a}e_{\b})e_{\g}-e_{\b}(e_{\g}e_{\a}),
\end{array}$$
which proves our claim. \qed

\medskip
This means that the model that we have found for $\fe_7$ really has 
a very close connection with the octonions. We can reformulate our
discussion as follows. 

\begin{theo}
The exceptional complex Lie algebra $\fe_7$ has a natural structure
of an $\OO$-graded algebra, given in terms of points $i$ 
and lines $\ell$ on the Fano plane by
$$\fe_7=\times_{i=1}^7\fsl(A_i)e_0\;\op\bigoplus_{\ell}
(\otimes_{i\notin\ell}A_i) \, e_{\ell}.$$
Moreover, its minimal representation decomposes as 
$$V=\bigoplus_{\ell}(\otimes_{i\in\ell}A_i) \, e_{\ell}.$$
\end{theo}

Note that there is a quaternionic analogue of this construction,
where instead of the Fano plane we consider only one of its lines. 
This means that we glue three copies of $\fso_8$ along the product 
of four copies of $\fsl_2$. The resulting algebra is $\fso_{12}$. 

We conclude that the Lie algebra $\fe_7$ is the support of a 
finite plane projective geometry whose points represent $7$ 
copies of $\fso_8$, and whose lines represent  $7$ copies of $\fso_{12}$. 

\subsection*{A sign problem and the isomorphism $PSL(3,\FF_2)\simeq
PSL(2,\FF_7)$}
We have just 
seen that the octonionic multiplication table gives a solution 
to the problem of finding a set of constants $\theta_{\a,\b}$ satisfying
the skewsymmetry condition and the $56$ quadratic relations associated
to the $28$ triangles in the Fano plane. What are the other solutions 
such that $\theta_{\a,\b}=\pm 1$? 

\begin{prop} There exist exactly sixteen such solutions, falling into 
two $PSL(3,\FF_2)$-orbits. Each orbit can be identified, as a 
$PSL(2,\FF_7)$-set, with a copy of the projective line $\FF_7\PP^1$. 
\end{prop}

\proof We proceed as follows. We first check that on a line the 
orientations must be coherent in the following sense: put an arrow
from $\a$ to $\b$ if $\theta_{\a\b}=+1$. Then the three arrows on 
a line, if we draw it as a circle, must go in the same direction. In 
particular there are only two possible choices of signs on a line, 
one for each cyclic orientation. We can switch from one to the other
by changing one of the basis vectors in its opposite. Moreover, 
the possible solutions to our problem can now be interpreted 
as a coherent orientation of the seven lines in the plane. 

Now we choose a triangle in the Fano plane. We have eight possible 
orientations for the three sides. We observe that once we choose one, 
the orientation of the line joining the three middle points of the 
sides of the triangle is fixed, and that there are only two
possibilities for the three remaining lines, those passing 
through the center of the triangle. Moreover we pass from one
to the other by changing the sign of the basis vector corresponding
to the center. 

Finally we check that once we fix a coherent orientation, we can
transform it by $PSL_3(\FF_2)$ to an arbitrarily chosen orientation 
on the triangle of reference. This implies that we have $16$ possible
orientations splitting in at most two orbits. But there cannot be a 
single orbit since $16$ does not divide the order of $PSL_3(\FF_2)$.
To identify the two orbits with a projective line over $\FF_7$, 
there just remains to observe that $PSL_3(\FF_2)$, up to conjugation, has a
unique subgroup of index $8$ (see \cite{atlas}). \qed

\medskip
This suggests an interpretation of the isomorphism 
between $PSL(3,\FF_2)$ and $PSL(2,\FF_7)$. We obtained the following
one which we could not find in the litterature. It is conveniently 
expressed in terms of cubes in a projective line, by which we simply
mean a graph with eight vertices, which is topologically the 
incidence graph of a cube. 

\medskip\noindent 
{\it Definition}. A cube in a projective line $\PP^1$
is {\it harmonic} if each of its faces $(wxyz)$ is harmonic, 
that is, the opposite vertices $(wy)$ and $(xz)$ are in 
harmonic position. 

\begin{prop}
There exist fourteen harmonic cubes in $\FF_7\PP^1$, made of $42$
harmonic faces. They split uniquely into two $PSL(2,\FF_7)$-orbits 
in such a way that each harmonic face belongs to exactly one 
cube of each family.  
\end{prop}

\medskip
\begin{center}
\setlength{\unitlength}{3.5mm}
\begin{picture}(20,42)(-1,-1)
\multiput(0,0)(0,3){2}{\line(1,0){3}} 
\multiput(0,0)(3,0){2}{\line(0,1){3}} 
\multiput(2,1)(0,3){2}{\line(1,0){3}} 
\multiput(2,1)(3,0){2}{\line(0,1){3}} 
\multiput(0,0)(0,3){2}{\line(2,1){2}} 
\multiput(3,0)(0,3){2}{\line(2,1){2}} 
\put(-.8,-.7){5}\put(3.2,-.7){3}\put(5.3,.5){2}\put(1.3,1.2){6}
\put(-.8,3.1){0}\put(3.3,2.5){$\infty$}\put(5.2,4.2){1}\put(1.3,4.2){4}
\put(-4,2){$(1)$}

\multiput(16,0)(0,3){2}{\line(1,0){3}} 
\multiput(16,0)(3,0){2}{\line(0,1){3}} 
\multiput(18,1)(0,3){2}{\line(1,0){3}} 
\multiput(18,1)(3,0){2}{\line(0,1){3}} 
\multiput(16,0)(0,3){2}{\line(2,1){2}} 
\multiput(19,0)(0,3){2}{\line(2,1){2}} 
\put(15.2,-.7){6}\put(19.2,-.7){5}\put(21.3,.5){3}\put(17.3,1.2){2}
\put(15.2,3.1){0}\put(19.3,2.5){$\infty$}\put(21.2,4.2){1}\put(17.3,4.2){4}
\put(10,2){$(123)$}

\multiput(0,6)(0,3){2}{\line(1,0){3}} 
\multiput(0,6)(3,0){2}{\line(0,1){3}} 
\multiput(2,7)(0,3){2}{\line(1,0){3}} 
\multiput(2,7)(3,0){2}{\line(0,1){3}} 
\multiput(0,6)(0,3){2}{\line(2,1){2}} 
\multiput(3,6)(0,3){2}{\line(2,1){2}} 
\put(-.8,5.3){2}\put(3.2,5.3){5}\put(5.3,6.5){$\infty$}\put(1.3,7.2){6}
\put(-.8,9.1){4}\put(3.3,8.3){3}\put(5.2,10.2){1}\put(1.3,10.2){0}
\put(-4,8){$(2)$}

\multiput(16,6)(0,3){2}{\line(1,0){3}} 
\multiput(16,6)(3,0){2}{\line(0,1){3}} 
\multiput(18,7)(0,3){2}{\line(1,0){3}} 
\multiput(18,7)(3,0){2}{\line(0,1){3}} 
\multiput(16,6)(0,3){2}{\line(2,1){2}} 
\multiput(19,6)(0,3){2}{\line(2,1){2}} 
\put(15.2,5.3){5}\put(19.2,5.3){6}\put(21.3,6.5){1}\put(17.3,7.2){2}
\put(15.2,9.1){0}\put(19.3,8.3){4}\put(21.2,10.2){$\infty$}\put(17.3,10.2){3}
\put(10,8){$(147)$}

\multiput(0,12)(0,3){2}{\line(1,0){3}} 
\multiput(0,12)(3,0){2}{\line(0,1){3}} 
\multiput(2,13)(0,3){2}{\line(1,0){3}} 
\multiput(2,13)(3,0){2}{\line(0,1){3}} 
\multiput(0,12)(0,3){2}{\line(2,1){2}} 
\multiput(3,12)(0,3){2}{\line(2,1){2}} 
\put(-.8,11.3){1}\put(3.2,11.3){4}\put(5.3,12.5){$\infty$}\put(1.3,13.2){5}
\put(-.8,15.1){3}\put(3.3,14.3){2}\put(5.2,16.2){0}\put(1.3,16.2){6}
\put(-4,14){$(3)$}

\multiput(16,12)(0,3){2}{\line(1,0){3}} 
\multiput(16,12)(3,0){2}{\line(0,1){3}} 
\multiput(18,13)(0,3){2}{\line(1,0){3}} 
\multiput(18,13)(3,0){2}{\line(0,1){3}} 
\multiput(16,12)(0,3){2}{\line(2,1){2}} 
\multiput(19,12)(0,3){2}{\line(2,1){2}} 
\put(15.2,11.3){6}\put(19.2,11.3){2}\put(21.3,12.5){$\infty$}\put(17.3,13.2){3}
\put(15.2,15.1){4}\put(19.3,14.3){1}\put(21.2,16.2){0}\put(17.3,16.2){5}
\put(10,14){$(156)$}

\multiput(0,18)(0,3){2}{\line(1,0){3}} 
\multiput(0,18)(3,0){2}{\line(0,1){3}} 
\multiput(2,19)(0,3){2}{\line(1,0){3}} 
\multiput(2,19)(3,0){2}{\line(0,1){3}} 
\multiput(0,18)(0,3){2}{\line(2,1){2}} 
\multiput(3,18)(0,3){2}{\line(2,1){2}} 
\put(-.8,17.3){5}\put(3.2,17.3){4}\put(5.3,18.5){$\infty$}\put(1.3,19.2){6}
\put(-.8,21.1){2}\put(3.3,20.3){0}\put(5.2,22.2){3}\put(1.3,22.2){1}
\put(-4,20){$(4)$}

\multiput(16,18)(0,3){2}{\line(1,0){3}} 
\multiput(16,18)(3,0){2}{\line(0,1){3}} 
\multiput(18,19)(0,3){2}{\line(1,0){3}} 
\multiput(18,19)(3,0){2}{\line(0,1){3}} 
\multiput(16,18)(0,3){2}{\line(2,1){2}} 
\multiput(19,18)(0,3){2}{\line(2,1){2}} 
\put(15.2,17.3){4}\put(19.2,17.3){5}\put(21.3,18.5){$\infty$}\put(17.3,19.2){3}
\put(15.2,21.1){0}\put(19.3,20.3){2}\put(21.2,22.2){6}\put(17.3,22.2){1}
\put(10,20){$(246)$}

\multiput(0,24)(0,3){2}{\line(1,0){3}} 
\multiput(0,24)(3,0){2}{\line(0,1){3}} 
\multiput(2,25)(0,3){2}{\line(1,0){3}} 
\multiput(2,25)(3,0){2}{\line(0,1){3}} 
\multiput(0,24)(0,3){2}{\line(2,1){2}} 
\multiput(3,24)(0,3){2}{\line(2,1){2}} 
\put(-.8,23.3){4}\put(3.2,23.3){2}\put(5.3,24.5){$\infty$}\put(1.3,25.2){6}
\put(-.8,27.1){5}\put(3.3,26.3){1}\put(5.2,28.2){0}\put(1.3,28.2){3}
\put(-4,26){$(5)$}

\multiput(16,24)(0,3){2}{\line(1,0){3}} 
\multiput(16,24)(3,0){2}{\line(0,1){3}} 
\multiput(18,25)(0,3){2}{\line(1,0){3}} 
\multiput(18,25)(3,0){2}{\line(0,1){3}} 
\multiput(16,24)(0,3){2}{\line(2,1){2}} 
\multiput(19,24)(0,3){2}{\line(2,1){2}} 
\put(15.2,23.3){6}\put(19.2,23.3){1}\put(21.3,24.5){4}\put(17.3,25.2){3}
\put(15.2,27.1){0}\put(19.3,26.3){5}\put(21.2,28.2){$\infty$}\put(17.3,28.2){2}
\put(10,26){$(367)$}

\multiput(0,30)(0,3){2}{\line(1,0){3}} 
\multiput(0,30)(3,0){2}{\line(0,1){3}} 
\multiput(2,31)(0,3){2}{\line(1,0){3}} 
\multiput(2,31)(3,0){2}{\line(0,1){3}} 
\multiput(0,30)(0,3){2}{\line(2,1){2}} 
\multiput(3,30)(0,3){2}{\line(2,1){2}} 
\put(-.8,29.3){4}\put(3.2,29.3){5}\put(5.3,30.5){$\infty$}\put(1.3,31.2){3}
\put(-.8,33.1){1}\put(3.3,32.3){0}\put(5.2,34.2){2}\put(1.3,34.2){6}
\put(-4,32){$(6)$}

\multiput(16,30)(0,3){2}{\line(1,0){3}} 
\multiput(16,30)(3,0){2}{\line(0,1){3}} 
\multiput(18,31)(0,3){2}{\line(1,0){3}} 
\multiput(18,31)(3,0){2}{\line(0,1){3}} 
\multiput(16,30)(0,3){2}{\line(2,1){2}} 
\multiput(19,30)(0,3){2}{\line(2,1){2}} 
\put(15.2,29.3){4}\put(19.2,29.3){2}\put(21.3,31.5){$\infty$}\put(17.3,31.2){6}
\put(15.2,33.1){3}\put(19.3,32.3){5}\put(21.2,34.2){1}\put(17.3,34.2){0}
\put(10,32){$(257)$}

\multiput(0,36)(0,3){2}{\line(1,0){3}} 
\multiput(0,36)(3,0){2}{\line(0,1){3}} 
\multiput(2,37)(0,3){2}{\line(1,0){3}} 
\multiput(2,37)(3,0){2}{\line(0,1){3}} 
\multiput(0,36)(0,3){2}{\line(2,1){2}} 
\multiput(3,36)(0,3){2}{\line(2,1){2}} 
\put(-.8,35.3){3}\put(3.2,35.3){4}\put(5.3,36.5){$\infty$}\put(1.3,37.2){2}
\put(-.8,39.1){0}\put(3.3,38.3){6}\put(5.2,40.2){1}\put(1.3,40.2){5}
\put(-4,38){$(7)$}

\multiput(16,36)(0,3){2}{\line(1,0){3}} 
\multiput(16,36)(3,0){2}{\line(0,1){3}} 
\multiput(18,37)(0,3){2}{\line(1,0){3}} 
\multiput(18,37)(3,0){2}{\line(0,1){3}} 
\multiput(16,36)(0,3){2}{\line(2,1){2}} 
\multiput(19,36)(0,3){2}{\line(2,1){2}} 
\put(15.2,35.3){5}\put(19.2,35.3){4}\put(21.3,36.5){$\infty$}\put(17.3,37.2){6}
\put(15.2,39.1){1}\put(19.3,38.3){2}\put(21.2,40.2){0}\put(17.3,40.2){3}
\put(10,38){$(345)$}

\end{picture} 
\end{center}
\centerline{{\it Figure 4. The fourteen harmonic cubes in $\FF_7\PP^1$}}

\medskip
We thus get an interpretation of these two sets of harmonic
cubes in $\FF_7\PP^1$ as points and lines in a Fano plane. 
We arbitrarily distinguish these two sets by calling them 
$p$-cubes and $\ell$-cubes, respectively. (Note an interesting 
analogy with the construction of the projective space over
$\FF_2$ given in \cite{polster} from the Fano plane: there 
exists $30$ unequivalent labelings of the seven vertices up to
the $PSL(3,\FF_2)$ action, and they are split into two families 
of $15$ labelings by the property that each line appears only once 
in each family.) 

The incidence relations can be defined as follows:
\begin{itemize}
\item Given a $p$-cube (respectively $\ell$-cube), there exist
exactly three $\ell$-cubes (respectively $p$-cubes) sharing a pair
of opposite faces with it. 
\item Given two $p$-cubes (respectively $\ell$-cubes), there exists 
a unique $p$-cube (respectively $\ell$-cube) such that the three cubes 
can be split each into two pairs of opposite edges forming squares
with the same fourtuples of vertices.  
\end{itemize}

Note that each pair $(ij)$ is the diagonal of exactly one $p$-cube 
and one $\ell$-cube. 

\smallskip
A nice feature of the correspondence is that each pair of points 
in $\FF_7\PP^1$ defines an edge of exactly three cubes in each family, 
corresponding to the three vertices and to the three edges of a 
triangle in the Fano plane. Therefore:

\begin{prop}
There is an equivariant correspondence between the $28$ triangles
in the Fano plane, and the $28$ pairs of points in the projective
line over $\FF_7$.
\end{prop}

Explicitely, this correspondence is as follows, where the triples in 
boldface are the vertices of a triangle in the Fano plane:

$$\begin{array}{llllllll}
01 & {\bf 256}\hspace*{5mm} & 12 & {\bf 145}\hspace*{5mm} & 24 
& {\bf 235}\hspace*{5mm} & 3\infty & {\bf 146} \\
02 & {\bf 346} & 13 & {\bf 234} & 25
 & {\bf 247} & 45 & {\bf 456} \\
03 & {\bf 457} & 14 & {\bf 136} & 26
 & {\bf 126} & 46 & {\bf 157} \\
04 & {\bf 124} & 15 & {\bf 357} & 2\infty 
& {\bf 567} & 4\infty & {\bf 347}  \\
05 & {\bf 167} & 16 & {\bf 467} & 34
 & {\bf 267} & 56 & {\bf 134} \\
06 & {\bf 237} & 1\infty & {\bf 127} & 35
 & {\bf 125} & 5\infty & {\bf 236} \\
0\infty & {\bf 135} & 23 & {\bf 137} & 36
 & {\bf 356} & 6\infty & {\bf 245}
\end{array}$$

\medskip
What about the orientations of the Fano planes that we were
interested in? We can associate such an orientation, in a 
$PSL(2,\FF_7)$-equivariant way, to each point $p\in\FF_7\PP^1$
as follows. For each point $x$ in the Fano plane, consider the 
point $q_x\in\FF_7\PP^1$ such that $pq_x$ is a diagonal of the 
harmonic cube corresponding to $x$. (This defines a bijection 
between $\FF_2\PP^2$ and $\FF_7\PP^1-\{p\}$.) Then the line 
$\ell=(xyz)$ will we positively oriented if the cross-ratio
$$(pq_xq_yq_z)=3.$$
On can easily check that in $\FF_7$, this condition is 
invariant under a cyclic permutation of $(xyz)$, so that this
definition really makes sense! This makes explicit the identification 
that we obtained between the projective line over $\FF_7$ and
half of the coherent orientations of the Fano plane. Of course 
we obtain the other half by reversing all the arrows. 

\medskip
Our correspondence has the property to transform certain
special configurations of bitangents into nice sets of triangles. 
We mention three instances of that. 

\smallskip
{\bf 1}. When we index the bitangents by pairs of points in a set with eight
elements, we give a special role to the eight Aronhold sets 
$A_i$ defined as the seven bitangents $(ij)$, $j\ne i$. Indeed,
this follows from the discussion of Example 4 above. We thus get
eight sets of seven triangles $T_i$ in the Fano plane (corresponding to
the eight points in $\FF_7\PP^1$) with the 
following properties: given a line $\ell$ and two points on it, 
there is a unique triangle in $T_i$ they are vertices of which;
in particular, $T_i$ is a copy of the Steiner triple system
$S(2,3,7)$; from the three pairs of points in $\ell$ we thus deduce three 
triangles in $T_i$; the three vertices of these triangles which
do not belong to $\ell$ are the vertices of a fourth triangle in 
$T_i$; this defines a natural bijection between the lines of the Fano
plane and the triangles in $T_i$.  

$$\begin{array}{llllllll}
 & {\bf 123} & {\bf 174} & {\bf 156} & {\bf 246} & {\bf 257}
 & {\bf 345} & {\bf 376} \\
{\bf T_0} & 475 & 265 & 273 & 135 & 364 & 167 & 142 \\ 
{\bf T_1} & 467 & 265 & 234 & 375 & 163 & 172 & 154 \\ 
{\bf T_2} & 576 & 253 & 247 & 173 & 364 & 126 & 154 \\ 
{\bf T_3} & 475 & 356 & 234 & 173 & 146 & 276 & 125 \\ 
{\bf T_4} & 456 & 253 & 374 & 157 & 163 & 276 & 142 \\ 
{\bf T_5} & 456 & 236 & 247 & 375 & 134 & 167 & 125 \\ 
{\bf T_6} & 467 & 356 & 273 & 157 & 134 & 126 & 245 \\ 
{\bf T_{\infty}} & 576 & 236 & 374 & 135 & 465 & 172 & 245  
\end{array}$$

\medskip
This gives a remarkable configuration of 
eight Steiner triple systems formed on $28$ triangles
in such a way that each of them appears exactly twice.

Moreover, each point and each line in the Fano plane belongs to 
exactly three of the seven triangles in each system. And 
the centers of the 
seven triangles are the seven points of the plane.

\medskip {\bf 2}. We have seen in Example 5 
that syzygetic tetrads of bitangents are defined by fourtuples 
of weights in the fundamental representation of $\fe_7$ summing
to zero, but not in two opposite pairs. There are two types 
of such tetrads in Hesse's notations: $105$ are permutations 
of $(01)(23)(45)(6\infty)$ and $210$ are permutations of 
$(01)(23)(02)(13)$. Each of these tetrads defines a special 
configuration of four triangles, a typical one being   

\begin{center}
\setlength{\unitlength}{4mm}
\begin{picture}(10,8)(-3,-1)
\put(0,0){\line(1,0){6}}
\put(0,0){\line(1,2){3}}
\put(6,0){\line(-1,2){3}}
\put(1.5,3){\line(1,0){3}}
\put(1.5,3){\line(1,-2){1.5}}
\put(1.5,3){\line(3,-2){1.5}}
\put(3,0){\line(0,1){2}}
\put(4.5,3){\line(-1,-2){1.5}}
\put(4.5,3){\line(-3,-2){1.5}}
\end{picture}
\end{center}

The four triangles in this picture are the great exterior
triangle, plus the three triangles having the center of the 
picture for vertex, plus two others on the middle points of 
two of the sides of the first one.  

\medskip {\bf 3}. A triple of points on the projective line over
$\FF_7$ defines three pairs, hence three triangles in the Fano plane. 
The following statement leads us back to the orientation 
problem which was the starting point of this long digression. 

\begin{prop}
There is an induced equivariant correspondence between triples of points 
on the projective line over $\FF_7$, and oriented triangles in the 
Fano plane.
\end{prop}

This goes as follows. Consider an oriented triangle in $\FF_2\PP^2$. 
For each vertex, consider the middle point on the opposite side, 
and then go to the next vertex following the orientation. We thus get
three (non oriented) triangles, which can be checked to be in
correspondence with three pairs $(ab),(bc),(ac)$ of a unique triple
$(abc)$ of points on $\FF_7\PP^1$. 

Given a pair $(pq)$ of points in $\FF_7\PP^1$, it defines a triangle in the
Fano plane, hence two oriented triangles, hence two triples of points
in the projective line. How can we obtain them directly? Simply by 
considering the unique harmonic cube in one of our two families 
having $(pq)$ for diagonal. Then the three points on an edge of 
this cube passing to $p$ (respectively $q$) give the two triples. 

\subsection*{Another setting for the bitangents}
The natural inclusion $SL(3,\FF_2)\subset Sp(6,\FF_2)=W(E_7)^+$
suggests to encode the $28$ bitangents and their symmetry group,
directly in the geometry of the Fano plane. This can indeed be
done in a very natural way, once we have identified the 
bitangents with the $28$ triangles in the Fano plane. Recall that
the group of the bitangents is generated by the transpositions
$s_{ij}$ and the bifid tranformations $s_{(pqrs|xyzt)}$. We have
checked that they have a simple geometric interpretation in 
terms of triangles. 

\smallskip\noindent {\it Transpositions}.
Let $T$ be the triangle associated to the pair $(ij)$. There is a unique 
point in the Fano plane which does not belong to a side of $T$, and we
call this point the center of the triangle. Up to the action of 
$SL(3,\FF_2)$ we can draw the Fano plane in such a way that the 
exterior of the picture is precisely $T$. Then the involution $\s_{T}$
on the set of triangles exchanges triangles as shown in the following
picture: 

\begin{center}
\setlength{\unitlength}{4mm}
\begin{picture}(20,8)(-1,-1)
\put(0,0){\line(1,0){6}}
\put(0,0){\line(1,2){3}}
\put(6,0){\line(-1,2){3}}
\put(12,0){\line(1,0){6}}
\put(12,0){\line(1,2){3}}
\put(18,0){\line(-1,2){3}}
\put(8.6,3.6){$\sigma_T$}  

\thicklines
\put(0,0){\line(1,0){3}}
\put(3,0){\line(0,1){2}}
\put(0,0){\line(3,2){3}}
\put(15,0){\line(0,1){6}}
\put(13.5,3){\line(1,2){1.5}}
\put(13.5,3){\line(1,-2){1.5}}
\put(7.5,3){\vector(1,0){3}}
\put(10.5,3){\vector(-1,0){3}}

\end{picture}
\end{center}

Otherwise said, a triangle whose vertices are a vertex $v$ 
of $T$, the center $c$ of $T$ and the middle point $p$ of a side 
of $T$, is mapped to the triangle whose vertices are $p$, its 
symmetric point $w$ with respect to $c$, and the middle point of the 
side $vw$ -- and conversely, while the other triangles remain
unchanged. 

\smallskip\noindent {\it Bifid transformations}.
One can check that the $35$ bifid transformations, when we 
interprete them as operations on the triangles, split into 
three types which are naturally associated to the seven points 
$p$, the seven lines $\ell$, and the $21$ pairs of incident points 
and lines $p\in\ell$ in the Fano plane. We get the following
transformations. 

Associated to a point $p$ is a transformation $\sigma_p$ who
takes a triangle with a side whose middle point is $p$, and 
changes the opposite vertex to the symmetric point with respect to $p$.  

\begin{center}
\setlength{\unitlength}{4mm}
\begin{picture}(20,8)(-1,-1)
\put(0,0){\line(1,0){6}}
\put(0,0){\line(1,2){3}}
\put(6,0){\line(-1,2){3}}
\put(12,0){\line(1,0){6}}
\put(12,0){\line(1,2){3}}
\put(18,0){\line(-1,2){3}}
\put(8.6,3.6){$\sigma_p$}  
\put(-.8,-.5){$p$}
\put(11.2,-.5){$p$}

\thicklines
\put(3,0){\line(1,0){3}}
\put(3,0){\line(1,2){1.5}}
\put(6,0){\line(-1,2){1.5}}
\put(15,0){\line(1,0){3}}
\put(15,0){\line(0,1){2}}
\put(15,2){\line(3,-2){3}}
\put(7.5,3){\vector(1,0){3}}
\put(10.5,3){\vector(-1,0){3}}

\end{picture}
\end{center}

Associated to a line $\ell$ is a transformation $\sigma_{\ell}$ 
who takes a triangle with a unique vertex $v$ on $\ell$, and 
changes the two other vertices to the symmetric points with 
respect to $v$.  

\begin{center}
\setlength{\unitlength}{4mm}
\begin{picture}(20,8)(-1,-1)
\put(0,0){\line(1,0){6}}
\put(0,0){\line(1,2){3}}
\put(6,0){\line(-1,2){3}}
\put(12,0){\line(1,0){6}}
\put(12,0){\line(1,2){3}}
\put(18,0){\line(-1,2){3}}
\put(8.6,3.6){$\sigma_{\ell}$}  
\put(2.8,-.9){$\ell$}
\put(14.8,-.9){$\ell$}

\thicklines
\put(1.5,3){\line(3,-2){4.5}}
\put(1.5,3){\line(1,2){1.5}}
\put(3,6){\line(1,-2){3}}
\put(15,2){\line(3,-2){3}}
\put(15,2){\line(3,2){1.5}}
\put(16.5,3){\line(1,-2){1.5}}
\put(7.5,3){\vector(1,0){3}}
\put(10.5,3){\vector(-1,0){3}}

\end{picture}
\end{center}

Note that these two types of transformations are exchanged by the 
polarity transformation on the set of triangles, which exchanges
vertices and sides in the Fano plane with sides and vertices in the dual
Fano plane. 

\smallskip 
Finally, associated to a pair $p\in\ell$ is a transformation 
$\sigma_{p,\ell}$ who takes a triangle with a unique vertex 
$v\ne p$ on $\ell$, whose opposite side $s$ has $p$ for middle point, 
to the triangle with vertex the symmetric point of $v$ with respect 
to $p$, and opposite side the symmetric of $s$ with respect to $\ell$. 

\begin{center}
\setlength{\unitlength}{4mm}
\begin{picture}(20,8)(-1,-1)
\put(0,0){\line(1,0){6}}
\put(0,0){\line(1,2){3}}
\put(6,0){\line(-1,2){3}}
\put(12,0){\line(1,0){6}}
\put(12,0){\line(1,2){3}}
\put(18,0){\line(-1,2){3}}
\put(8.4,3.6){$\sigma_{p,\ell}$}  
\put(2.8,-.9){$\ell$}
\put(14.8,-.9){$\ell$}
\put(-.8,-.5){$p$}
\put(11.2,-.5){$p$}

\thicklines
\put(3,0){\line(0,1){6}}
\put(1.5,3){\line(1,2){1.5}}
\put(1.5,3){\line(1,-2){1.5}}
\put(15,2){\line(3,-2){3}}
\put(15,2){\line(3,2){1.5}}
\put(16.5,3){\line(1,-2){1.5}}
\put(7.5,3){\vector(1,0){3}}
\put(10.5,3){\vector(-1,0){3}}

\end{picture}
\end{center}

The group of the bitangents is isomorphic with the group of
permutations of the triangles generated by the elementary 
tranformations $\s_T$, $\s_p$, $\s_{\ell}$, $\s_{p,\ell}$. 
It contains $SL(3,\FF_2)$ as the group of collinations acting
on the triangles. This makes clear the natural inclusions
$$SL(3,\FF_2)\subset\cS_8\subset Sp(6,\FF_2)=W(E_7)^+.$$

\medskip

\section{Del Pezzo surfaces of degree one}

The bicanonical model of a Del Pezzo surface of degree one is 
the double cover of 
a quadratic cone, branched over a canonical space curve 
of genus $4$ and degree $6$ given by the complete
intersection of the cone with a unique cubic surface
(\cite{dem} V.1). 
The $240$ lines on the Del Pezzo surface arise in pairs 
from the $120$ tritangent planes to the canonical curve, 
which can be identified with its odd theta 
characteristics. Moreover, the fact 
that the unique quadric containing this curve is a cone
distinguishes one of the $136$ even theta characteristic
by the property that it vanishes at the vertex of that cone. 

The automorphism group of the $240$ lines is the Weyl group 
$W(E_8)$ of the root system of type $E_8$. Its order is  
$696\, 729\, 600 = 128\times 27\times 5 \times 8! =2^{13}3^{5}5^27$.
The automorphism group of the $120$ tritangent planes is the 
quotient by the normal subgroup $\{\pm 1\}$ and can be identified 
with the orthogonal group $O(8,\FF_2)^+$ which preserves the quadratic 
form given by half the natural norm on the root lattice mod 2
(\cite{bou}, Exercice 1 of section 4, p. 228). Among the $2^8=256$
points in $\FF_2^8$, those with norm one can be identified with the 
odd theta-characteristics,  and those with norm zero with the 
even theta-characteristics, including the special one which 
identifies with the origin. 

If we consider the action of the adjoint group
$E_8$ on the projectivized adjoint representation $\PP\fe_8$,
the $240$ root spaces can be interpreted as a kind of finite
skeleton of the closed orbit, the adjoint variety of $E_8$.
This variety parametrizes the so-called symplecta in 
Freudenthal's metasymplectic geometry (see \cite{LMfreud}). 

\smallskip
The affine Dynkin diagram of $E_8$ is 

\begin{center}
\setlength{\unitlength}{3mm}
\begin{picture}(30,5)(-8,-3)

\multiput(0,0)(2,0){8}{$\circ$}
\multiput(0.5,.4)(2,0){7}{\line(1,0){1.6}} 
\put(4,-2){$\circ$}
\put(4.3,-1.4){\line(0,1){1.5}}

\end{picture} 
\end{center}

Relevant configurations will be provided by the simplest 
markings.

\medskip\noindent {\it Example 1}. We mark the node next to the
rightmost extremal node. Then $\fh=\fsl_2\times \fe_7$ and $W$ has index
$120$ in $W(E_8)$.  

\begin{center}
\setlength{\unitlength}{3mm}
\begin{picture}(30,5)(-8,-3)

\multiput(0,0)(2,0){8}{$\circ$}
\multiput(0.5,.4)(2,0){7}{\line(1,0){1.6}} 
\put(4,-2){$\circ$}
\put(12,0){$\bullet$}
\put(4.3,-1.4){\line(0,1){1.5}}

\end{picture} 
\end{center}

The branching gives the binary model 
$$\fe_8=\fsl_2\times\fe_7\op A\otimes V,$$
where $V$ is again the minimal $56$-dimensional representation
of $\fe_7$. In particular, this associates to each root $\a$
(up to sign) of $\fe_8$ a complex $S_{\a}$ of $56$ tritangent planes. 
Two complexes $S_{\a}$ and $S_{\b}$ have two possible relative 
positions, distinguished by the fact that $(\a,\b)=0$ or not. 

In the latter case, exactly as for $\fe_7$ the two complexes are
{\it azygetic} and can be completed uniquely with a third 
complex $S_{\g}$, with $\g=\a\pm\b$, azygetic to both of them. 
There exists $1120$ such {\it azygetic triads of complexes}. 
This leads to the ternary model of $\fe_8$,
$$\fe_8=\fsl_3\times\fe_6\op (B\otimes J)\op (B\otimes J)^*.$$
The $120$ tritangent planes are partitioned into 
the triple $(\a,\b,\g)$, three sets $S_{\a\b}$, $S_{\a\g}$, $S_{\b\g}$
of cardinality $27$, and their complement of cardinality
$36$, with
$$S_{\a}=\{\b,\g\}\cup S_{\a\b}\cup S_{\a\g}.$$

In the former case, the two complexes are {\it syzygetic}. 
Since their common orthogonal subsystem is a root system of type 
$A_1\times A_1\times D_4$, we can define a {\it syzygetic
tetrad of complexes} to consist in four syzygetic complexes
orthogonal to a $D_4$-subsystem of the root system $E_8$. 
Note that we have three syzygetic tetrads for each 
$D_4$-subsystem, making a total of $9450$ such tetrads. 

A pair of syzygetic complexes can be completed uniquely into 
a syzygetic tetrad $(\a,\b,\g,\d)$. The other $116$ tritangent planes 
are then partitioned into a set $S_{\a\b\g\d}$ of cardinality $8$, 
six sets $S_{\a\b}$, $S_{\a\g}$, $S_{\a\d}$, $S_{\b\g}$, $S_{\b\d}$, 
$S_{\g\d}$ of cardinality $16$, and their complement of cardinality
$12$. Here $$S_{\a}=S_{\a\b\g\d}\cup S_{\a\b}\cup S_{\a\g}
\cup S_{\a\d}.$$

\medskip\noindent {\it Example 2}. We mark the leftmost node.
Then $\fh=\fspin_{16}$ and $W$ has index $135$ in $W(E_8)$. 

\begin{center}
\setlength{\unitlength}{3mm}
\begin{picture}(30,5)(-8,-3)

\multiput(0,0)(2,0){8}{$\circ$}
\multiput(0.5,.4)(2,0){7}{\line(1,0){1.6}} 
\put(4,-2){$\circ$}
\put(0,0){$\bullet$}
\put(4.3,-1.4){\line(0,1){1.5}}

\end{picture} 
\end{center}

The branching gives another very nice model, 
$$\fe_8=\fspin_{16}\op\Delta,$$
where $\Delta$ is a half-spin representation, of dimension
$128$. 

Geometrically, we get twelve dimensional quadrics in the 
adjoint variety. 

\medskip\noindent {\it Example 3}. We mark the lowest node.
Then $\fh=\fsl_9$ and $W=\cS_9$ has index $1920$. 

\begin{center}
\setlength{\unitlength}{3mm}
\begin{picture}(30,5)(-8,-3)

\multiput(0,0)(2,0){8}{$\circ$}
\multiput(0.5,.4)(2,0){7}{\line(1,0){1.6}} 
\put(4,-2){$\bullet$}
\put(4.3,-1.4){\line(0,1){1.5}}

\end{picture} 
\end{center}

The branching gives one of the prettiest models of $\fe_8$,
namely 
$$\fe_8=\fsl_9\op\Lambda^3U\op\Lambda^6U,$$
where $U$ denotes the natural nine dimensional representation. 
The $\fsl_9$ factor defines what Du Val calls an $\a_8$ configuration 
in the uniform polytope $4_{21}$.  (Note that there are $960$ such 
configurations rather than $1920$. This is because of the invariance 
by $-1$, the longest element in $W(E_8)$: indeed its restriction to
the root system $A_8$ of $\fsl_9$ does not belong to $W(A_8)$ but 
defines an order two {\it outer} automorphism). It can be interpreted 
as a special system of then even theta-characteristics (\cite{duval}, 
p. 51). One deduces a special set of $84$, and the complementary
set of $36$ tritangent planes: the odd theta-characteristics 
of the former set are obtained as sums of three of the distinguished 
even theta-characteristics (which is visible from the fact that they 
correspond to weights of a third wedge power), and the latter as sums 
of only two of them. 

Geometrically, we obtain copies of $\PP^7$ in the adjoint variety. 
From this point of view there is a close analogy with 
Aronhold sets of bitangents to the plane quartic curve.  

There is also a strong analogy with the decomposition of $\fe_7$
that we described in Example 2 above, from which we recovered Hesse's
notation for the bitangents to the plane quartic, and the so called
 bifid transformations. Here the roots of $\fe_8$, up to sign, are split 
in two types by the model we are discussing. The roots that are 
weights of the factor $\Lambda^3U$ are indexed by a triple $(ijk)$
of integers between $1$ and $9$. The roots or $\fsl_9$ are indexed,
up to sign, by a pair $(ij)$ that we can identify with the triple $(0ij)$. 
The tritangent planes  are then indexed by the  $120$  triples
of integers between $1$ and $10$, a notation first introduced by Pascal 
(see \cite{coxe}). Moreover, the group of the tritangent 
planes is generated by the symmetric group $\cS_{9}$ plus the 
symmetries associated to the triples $(ijk)$, which are again called
{\it bifid transformations} (but beware that it does not contain 
the symmetric group $\cS_{10}$). An easy computation shows that 
the bifid transformation associated to $(ijk)$ exchanges 
$(ij\ell)$ with $(0k\ell)$ and sends a triple $(abc)$
to $(pqr)$ if $ijkabcpqr$ is a permutation of $123456789$ -- while 
the other triples are fixed.

\smallskip There is a connection with the model for $\fe_6$
discussed in Example 4 of section 3, as we can see by splitting 
the nine dimensional representation $U$ into three supplementary 
spaces of dimension three. In fact the corresponding subroot 
system of type $A_2^3$ in $E_8$ is orthogonal to another $A_2$
subsystem, and we get another interesting model for $\fe_8$, 
which is graded over $\ZZ_3\times\ZZ_3$:
$$\begin{array}{ll}
\fe_8  = & \fsl(A_1)\times\fsl(A_2)\times\fsl(A_3)\times\fsl(A_4)\op \\
 & \hspace*{2cm}
\op A_{1234}\op A_{1^*234}\op A_{12^*34}\op A_{1234^*} \\
 & \hspace*{2cm}
\op A_{1234}^*\op A_{1^*234}^*\op A_{12^*34}^*\op A_{1234^*}^*,
\end{array}$$
where $A_1,A_2,A_3,A_4$ are three dimensional and we used the 
notation $A_{1^*234}=A_1^*\ot A_2\ot A_3\ot A_4$. But this grading
does not seem to be supported by any interesting geometry.

\medskip\noindent {\it Example 2, continued}.
The trialitarian model of $\fe_8$ arises when we extract from
$\fspin_{16}$ two orthogonal copies of $\fspin_{8}$. This can 
be seen from the affine Dynkin diagram of type $D_8$:

\begin{center}
\setlength{\unitlength}{3mm}
\begin{picture}(30,5)(-8,-1.5)

\multiput(2,0)(2,0){5}{$\circ$}
\multiput(2.5,.4)(2,0){4}{\line(1,0){1.6}} 
\put(.1,1.8){$\circ$}
\put(.1,-1.8){$\circ$}
\put(.6,2){\line(1,-1){1.5}}
\put(.6,-1.4){\line(1,1){1.5}}
\put(6,0){$\bullet$}
\put(11.9,1.8){$\circ$}
\put(11.9,-1.8){$\circ$}
\put(12,2){\line(-1,-1){1.5}}
\put(12,-1.4){\line(-1,1){1.5}}

\end{picture} 
\end{center}
Indeed, we get 
$$\fe_8=\fspin_{8}\times\fspin_{8}\op (\OO_1\op\OO'_1)
\op (\OO_2\op\OO'_2)\op (\OO_3\op\OO'_3),$$
where $\OO_1,\OO_2,\OO_3$ are the three eight-dimensional
irreducible representations of $\fspin_8$. This models splits the $120$
tritangent planes into two groups of $12$ and three groups 
of $32$. With which geometric characterization? 

Then we can split each  $\fspin_{8}$ using four copies of 
$\fsl_2$; we know that each of its eight-dimensional representation
has a nice decomposition, and putting them together we obtain
$$\fe_8=\times_{i=1}^8\fsl(A_i)\op\bigoplus_{(ijkl)\in I}
A_i\ot A_j\ot A_k\ot A_l,$$
where $I$ is the following set of $14$ quadruples:
$$\begin{array}{ccc}
1234 & \qquad\qquad & 5678 \\
1256 & & 3478 \\
1278 & & 3456 \\
1357 & & 2468 \\
1368 & & 2457 \\
1458 & & 2367 \\
1467 & & 2358 
\end{array}$$
As in the case of $\fe_7$ these fourteen quadruples form a 
{\it Steiner quadruple system} $S(3,4,8)$. In fact there exists 
a unique such system up to isomorphism. 

\medskip\noindent {\it Remark.} Such a decomposition is induced by the 
choice of a root subsystem of type $A_1^8$ inside the root system 
$E_8$. This is equivalent to what Du Val calls a $\b_8$ configuration 
in the polytope $4_{21}$. There is also a correspondence with 
the $135$ even theta-characteristics, or the $135$ norm one vectors
in the root lattice mod 2, which can be seen as follows. 
Consider such a vector, for example $\th=\o_1+\o_2+\o_3+\o_4+
\o_5+\o_6+\o_7+\o_8$ if we denote by $\pm\o_i$ the weights of $A_i$. 
Among the $127$ projective lines through $\th$, $63$ are contained 
in the quadric, $56$ are tangents and only $8$ are true bisecants; on each of 
these bisecants there is a unique point outside the quadric and we
obtain a set of eight points defining a subsystem of type  $A_1^8$.  

\medskip
Let  us consider our Steiner quadruple system is some detail. 
We first note that each quadruple $(ijkl)\in I$ defines a copy of 
$\fso_8\simeq \fsl(A_i)\times \fsl(A_j)\times \fsl(A_k)\times \fsl(A_l)
\op (A_i\ot A_j\ot A_k\ot A_l)$ inside $\fe_8$. We have thus 
constructed $\fe_8$ by gluing together fourteen copies of $\fso_8$
overlapping over eight copies of $\fsl_2$!

\smallskip
We can interprete the   $14$ quadruples in $I$ as {\it points} of a 
configuration whose {\it lines} are triples of type $(ijkl),
(klmn),(ijmn)$. A straightforward inspection shows that there
are exactly $28$ lines. Moreover, each line has three points
and each point belongs to $6$ lines. In other words, we have
obtained a $(14_6,28_3)$-configuration.  

This configuration has the following interpretation. Consider
$\FF_2\PP^3$, the three dimensional projective space over the
field with two elements. It has fifteen points. Choose one, say
$p_{\infty}$, and throw it away. Since $\FF_2\PP^3$ contains
$35$ lines, $7$ of which passing through $p_{\infty}$, we remain
with $14$ points and $28$ lines whose incidence configuration is 
the one we are interested in. 

In particular, note the following properties:
\begin{enumerate}
\item each point $p$ has an antipodal point $p^*$, 
the unique point to which it is not joined by a line --
in $\FF_2\PP^3$, this is the third point of the line 
$\overline{pp_{\infty}}$;
\item the $6$ lines passing through a point $p$ split
naturally into three pairs, in such a way that the four points
different from $p$ on each pair, belong to a pair of
lines passing through $p^*$; in $\FF_2\PP^3$, these three pairs 
are cut by the three planes containing the line $\overline{pp_{\infty}}$; 
\item the configuration is made of eight copies of the Fano plane, 
corresponding to the eight planes in $\FF_2\PP^3$ not containing
$p_{\infty}$; the other seven planes give sub-configurations
of type $(6_2,4_3)$ which are pointed Fano planes.
\end{enumerate}

But the most satisfactory way to understand our configuration is 
probably to see it as a {\it doubled Fano plane}. By this we mean
that we can associate to each point of the Fano plane a pair 
of antipodal quadruples, in such a way that the $28$ lines of 
the configuration correspond four by four to the lines of the 
Fano plane. In terms of Lie subalgebras of $\fe_8$, this 
associates to each point of the Fano plane a copy of $\fso_8
\times\fso_8$, and to each line a copy of $\fso_{16}$. Moreover, 
there exists four copies of $\fso_{12}$ inside the $\fso_{16}$
defined by a line, meeting each $\fso_8\times\fso_8$ corresponding
to one of its points, along one of the two $\fso_8$ factors. 

We have even more structure if we note that each integer $i$
between $1$ and $8$ determines a copy of $\fe_7$ inside $\fe_8$. 
Two such copies meet along one of the $\fso_{12}$ indexed by 
the $28$ lines. This has an interesting combinatorial interpretation. 
Consider the following $8\times 8$ array:
$$\begin{array}{cccccccc}
0 & 1 & 2 & 3 & 4 & 5 & 6 & 7 \\
1 & 0 & 3 & 2 & 5 & 4 & 7 & 6 \\
2 & 3 & 0 & 1 & 6 & 7 & 4 & 5 \\
3 & 2 & 1 & 0 & 7 & 6 & 5 & 4 \\
4 & 5 & 6 & 7 & 0 & 1 & 2 & 3 \\
5 & 4 & 7 & 6 & 1 & 0 & 3 & 2 \\
6 & 7 & 4 & 5 & 2 & 3 & 0 & 1 \\
7 & 6 & 5 & 4 & 3 & 2 & 1 & 0  
\end{array}$$
It is symmetric with respect to the main diagonal, 
and each number $j$ between $0$
and $7$ appears once and only once in each line and in each column. 
If we associate to $j$ the four pairs $(lc)$ of numbers indexing
the lines and columns of the boxes where $j$ appears, this 
means that we obtain a partition of $[1,8]$ into four pairs, 
for a total number of $28$ distinct pairs:
$$\begin{array}{l}
(12)(34)(56)(78) \\
(13)(24)(57)(68) \\
(14)(23)(58)(67) \\
(15)(26)(37)(48) \\
(16)(25)(38)(47) \\
(17)(28)(35)(46) \\
(18)(27)(36)(45)
\end{array}$$
(Note that this array can also be obtained from the diagonals of 
each family of seven harmonic cubes of Figure 5.)
These seven partitions index our seven lines in the following way:
to each $(ij)(kl)(pq)(rs)$ are associated the six quadruples obtained
by selecting two of the four pairs. The incidence is given by the 
following rule: two lines being given, they can always be indexed 
by partitions of type  $(ij)(kl)(pq)(rs)$ $(ik)(jl)(pr)(qs)$, and 
then their intersection point is indexed by the pair of antipodal
quadruples $(ijkl)$ and $(pqrs)$.  

We recapitulate the main conclusions of our discussion:

\begin{theo}
The exceptional complex Lie algebra $\fe_8$ has a natural 
structure of an $\OO$-graded Lie algebra, obtained by 
gluing seven copies of $\fso_8\times\fso_8$ indexed by 
the points of a Fano plane. 

This occurs in such a way that the
three copies indexed by three points of a same line are glued
together in a copy of $\fso_{16}$.  
\end{theo}

\noindent {\it First application}.
This construction of $\fe_8$ isolates eight tritangent planes 
corresponding to the roots of the eight copies of $\fsl_2$, the
remaining $112$ tritangent planes being split into fourteen
groups of eight. 

Define a {\it syzygetic tetrad} of tritangent planes by the 
property that their twelve tangency points are the twelve intersection
points of the sextic canonical curve with some quadric hypersurface. 
This means that they can be defined by four roots summing to zero. 
We claim that the eight tritangent planes associated  to each factor 
$A_i\ot A_j\ot A_k\ot A_l$ in our decomposition of $\fe_8$, can be split 
in two syzygetic tetrads in exactly six ways. Indeed, the corresponding
roots are the $\pm\epsilon_i\pm\epsilon_j\pm\epsilon_k\pm\epsilon_l$
and two syzygetic tetrads are given by the following two sign tables:

$$\begin{array}{ccccccccc}
+ & + & + & + &\qquad  + & -  & + & - & \\
+ & + & - & - &\qquad  + & -  & -  & + & \\
- & - & + & - &\qquad  - & + & + & + & \\
- & - & - & + &\qquad  - & + & - & - & 
\end{array}$$ 

The other five are obtained by permuting the three last columns of
each table. 

\medskip\noindent {\it Second application}. 
An {\it orthogonal decomposition} (OD) of a semisimple Lie 
algebra $\fg$ is a decomposition $\fg=\oplus_{i=0}^h\ft_i$
into a direct sum of Cartan subalgebras. Such an OD is 
{\it multiplicative} if for each pair $(i,j)$, there exists
an integer $k$ such that $[\ft_i,\ft_j]\subset\ft_k$.
A trivial example is that of $\fsl_2=\fsl(A)$, once we have chosen
a basis $(e,f)$ of $A$. If $X,Y,H$ is the associated canonical 
basis of $\fsl_2(A)$, a multiplicative OD (or MOD) is given be the 
three lines generated by $H, X+Y$ and $X-Y$.  

Multiplicative OD's have been used by Thompson to construct 
the finite sporadic simple group denoted $Th$ or $F_3$. 
Indeed, his construction relied on the existence of a 
multiplicative OD for $\fe_8$. 

In fact there exists, up to isomorphism, a unique multiplicative
OD of $\fe_8$ (\cite{kos}, Chapter 3). Our construction provides 
it for free! Indeed
we just need to notice that each component $A_i\ot A_j\ot 
A_k\ot A_l$ can be split into the direct sum of two Cartan
subalgebras $\ft_{ijkl}^{\pm}$ by choosing a basis 
$(e_i,f_i)$ of each $A_i$, and letting 
$$\ft_{ijkl}^{\pm}=\langle x_i\ot x_j\ot x_k\ot x_l
\pm y_i\ot y_j\ot y_k\ot y_l\rangle,$$
where $(x_i,y_i)$ is $(e_i,f_i)$ or $(f_i,e_i)$.
This gives $28$ Cartan subalgebras of $\fe_8$, and we add three 
others, say $\ft_0$, $\ft_+$ and $\ft_-$, 
by putting together the OD's of the 
$\fsl(A_i)$'s associated to the basis of the $A_i$'s that  
we have chosen. We obtain:

\begin{theo}
The direct sum decomposition
$$\fe_8=\ft_0\op\ft_+\op\ft_-\op \bigoplus_{(ijkl)\in I}
(\ft_{ijkl}^+\op\ft_{ijkl}^-)$$
is a multiplicative OD of $\fe_8$.
\end{theo}

\bigskip

\providecommand{\bysame}{\leavevmode\hbox to3em{\hrulefill}\thinspace}
{\sc Laurent Manivel},

Institut Fourier, UMR 5582 (UJF-CNRS),

BP 74, 38402 St Martin d'H\`eres Cedex, France.

E-mail : {\tt Laurent.Manivel@ujf-grenoble.fr}

\end{document}